\documentclass[12pt]{article}
\usepackage{amsmath,amsfonts,latexsym,amsthm,amssymb}
\newcommand{\D}{\mathbb D_t^{(\alpha )}}
\numberwithin{equation}{section}

\DeclareMathOperator{\const}{const}

\DeclareMathOperator{\I}{Im}
\DeclareMathOperator{\Res}{Res}
\begin{document}
\newtheorem{prop}{Proposition}
\newtheorem{lem}{Lemma}
\newtheorem*{teo}{Theorem}
\pagestyle{plain}
\title{Cauchy Problem for Fractional Diffusion Equations}
\author{Samuil D. Eidelman and Anatoly N. Kochubei\footnote{Partially
supported by CRDF under Grant UM1-2421-KV-02}
\\ \footnotesize Institute of Mathematics,\\
\footnotesize National Academy of Sciences of Ukraine,\\
\footnotesize Tereshchenkivska 3, Kiev, 01601 Ukraine}
\date{}
\maketitle
\vspace*{2cm}
Running head:\quad  ``Fractional Diffusion Equations''

\vspace*{2cm}
The address for correspondence: \\

\vspace*{1cm}
Anatoly N. Kochubei\\
Institute of Mathematics,\\
National Academy of Sciences of Ukraine,\\
Tereshchenkivska 3, Kiev, 01601 Ukraine. \\
E-mail: \ kochubei@i.com.ua
\newpage
\vspace*{8cm}
\begin{abstract}
We consider an evolution equation with the regularized fractional
derivative of an order $\alpha \in (0,1)$ with respect to the time
variable, and a uniformly elliptic operator with variable coefficients
acting in the spatial variables. Such equations describe diffusion
on inhomogeneous fractals. A fundamental solution of the Cauchy
problem is constructed and investigated.
\end{abstract}
\vspace{2cm}
\newpage
\section{INTRODUCTION}

\medskip
In this paper we consider equations of the form
\begin{equation}
\left( \D u\right) (t,x)-Bu(t,x)=f(t,x),\quad t\in (0,T],\ x\in
\mathbb R^n,
\end{equation}
where $\D$, $0<\alpha <1$, is a regularized fractional derivative
(the Caputo derivative), that is
\begin{equation}
\left( \D u\right) (t,x)=\frac{1}{\Gamma (1-\alpha )}\left[
\frac{\partial}{\partial t}\int\limits_0^t(t-\tau )^{-\alpha
}u(\tau ,x)\,d\tau -t^{-\alpha }u(0,x)\right] ,
\end{equation}
$$
B=\sum\limits_{i,j=1}^na_{ij}(x)\frac{\partial^2}{\partial
x_i\partial x_j}+\sum\limits_{j=1}^nb_j(x)\frac{\partial
}{\partial x_j}+c(x)
$$
is a uniformly elliptic second order differential operator with
bounded continuous real-valued coefficients.

A strong motivation for investigating such equations comes from
physics. Fractional diffusion equations describe anomalous
diffusion on fractals (physical objects of fractional dimension,
like some amorphous semiconductors or strongly porous materials;
see \cite{AL,MK} and references therein). In normal diffusion
(described by the heat equation or more general parabolic
equations) the mean square displacement of a diffusive particle
behaves like $\const \cdot t$ for $t\to \infty$. A typical
behavior for anomalous diffusion is $\const \cdot t^\alpha$, and
this was the reason to invoke the equation (1.1), usually with
$B=\Delta$, where this anomalous behavior is an easy mathematical
fact. For connections to statistical mechanics see also
\cite{Gor,Meer}.

It is natural from the physical point of view to consider a usual
Cauchy problem, with the initial condition
\begin{equation}
u(0,x)=u_0(x).
\end{equation}
This setting determines the necessity to use the regularized fractional
derivative (1.2). If, for example, one considers instead of (1.2)
the Riemann-Liouville fractional derivative defined similarly, but
without subtracting $t^{-\alpha }u(0,x)$, then the appropriate
initial data will be the limit value, as $t\to 0$, of the
fractional integral of a solution of the order $1-\alpha$, not the
limit value of the solution itself (see e.g. \cite{SKM,Kost}).

Note that on a smooth enough function $v(t)$ the regularized fractional
derivative $\D v$ can be written as $(\mathbb D^{(\alpha
)}v)(t)=\dfrac{1}{\Gamma (1-\alpha
)}\int\limits_0^t\dfrac{v'(t)}{(t-\tau )^\alpha }d\tau$. In the
physical literature the expression on the right is
used as the basic object for formulating
fractional diffusion equations. However, in order to proceed
rigorously, one has either to use (1.2), or to consider, as in
\cite{SW}, an equivalent integral equation instead of the Cauchy
problem (1.1), (1.3). On the other hand, the above expression for
$\mathbb D^{(\alpha )}$ on smooth functions is a special case of a
generalized fractional differential operator of
Dzhrbashyan-Nersessyan \cite{DN} who studied ordinary differential
equations with such operators. For example, they showed that the
solution of the Cauchy problem $\mathbb D^{(\alpha )}v-\lambda
v=0$, $v(0)=1$, has the form $v(t)=E_\alpha (\lambda t^\alpha )$
where $E_\alpha$ is the Mittag-Leffler function. See also
\cite{MR}.

The mathematical theory of fractional diffusion equations has made
only its first steps. An expression for the fundamental solution
of the Cauchy problem (1.1), (1.3) with $B=\Delta$ was found
independently by Schneider and Wyss \cite{SW} and Kochubei
\cite{K2}. It was also shown in \cite{SW} that the fundamental
solution is non-negative, which led later \cite{S1,S2,Kols,Yor} to
a probabilistic interpretation of the equation (1.1). In \cite{SW}
only initial functions $u_0\in \mathcal S(\mathbb R^n)$ were
considered. A more general situation was studied in \cite{K2}
where $u_0$ was permitted to be unbounded, with minimal smoothness
assumptions. There are also some results regarding
initial-boundary value problems (see \cite{W,SW}).

For general problems (1.1), (1.3), in \cite{K2} a uniqueness
theorem for bounded solutions, and an exact uniqueness theorem
(for the case $n=1$) for solutions with a possible exponential
growth were proved; see Sect. 2 below. There are also several
papers devoted to the Cauchy problem for abstract evolution
equations (1.1), in which $B$ is a closed operator on a Banach
space (\cite{K1,El,Bazh,Bazh1,BM}, and others).

An obvious analogy with the classical theory of parabolic partial
differential equations \cite{E,F,LSU} suggests a wide range of
problems for the general equations (1.1) which deserve to be
investigated. It is natural to begin with the construction and
investigation of a fundamental solution of the Cauchy problem (1.1),
(1.3). That is the aim of this paper.

More specifically, we will construct and study, under natural
assumptions upon the coefficients of $B$, a Green matrix for the problem
(1.1), (1.3), that is such a pair
$$
\{ Z(t,x;\xi ),Y(t,x;\xi
)\},\quad (t\in (0,T],\  x,\xi \in \mathbb R^n)
$$
that for any bounded continuous function $u_0$ (locally H\"older
continuous, if $n>1$) and any bounded
function $f$, jointly continuous in $(t,x)$ and locally H\"older continuous
in $x$, there exists a classical
solution of the problem (1.1), (1.3) of the form
\begin{equation}
u(t,x)=\int\limits_{\mathbb R^n} Z(t,x;\xi )u_0(\xi )\,d\xi
+\int\limits_0^td\lambda \int\limits_{\mathbb R^n} Y(t-\lambda,x,;
y)f(\lambda ,y)\,dy.
\end{equation}

As in \cite{K2}, we call a function $u(t,x)$ a classical solution
if:
\begin{description}
\item[(i)] $u(t,x)$ is twice continuously differentiable in $x$
for each $t>0$;
\item[(ii)] for each $x\in \mathbb R^n$ $u(t,x)$ is continuous in
$t$ on $[0,T]$, and its fractional integral
$$
\left( I_{0+}^{1-\alpha }u\right) (t,x)=\frac{1}{\Gamma (1-\alpha
)}\int\limits_0^t(t-\tau )^{-\alpha }u(\tau ,x)\,d\tau
$$
is continuously differentiable in $t$ for $t>0$.
\item[(iii)]
$u(t,x)$ satisfies (1.1) and (1.3).
\end{description}

The results of this paper (as well as those of the earlier one
\cite{K2}) demonstrate a number of interesting features of the
equation (1.1), which represent a peculiar union of properties
typical for second order parabolic differential equations (like a
kind of the maximum principle; see the proofs of Theorems 1,2 in
\cite{K2}) and general parabolic equations and systems. It is well
known \cite{E,F,LSU} that the latter are characterized by the
parabolic weight $2b$ ($b\ge 1$), so that the differentiation in $t$
has the same ``force'' as the differentiation of the order $2b$
with respect to the spatial variables. In other words, the
differentiation in spatial variables has the weight $1/2b$ with
respect to $\partial /\partial t$. For the equation (1.1) we have
$\alpha /2$ instead of $1/2b$, and this formal analogy goes
through all the constructions and estimates including the
description of uniqueness and correctness classes.

In contrast to classical parabolic equations, the fundamental
solution $Z$ of the Cauchy problem (1.1), (1.3) (even for
$B=\Delta$) has, if $n>1$, a singularity not only at $t=0$, but
also at the ``diagonal'' $x=\xi$ with respect to the spatial
variables. This causes serious complications for the
implementation of the classical Levi method \cite{E,F,LSU}. The
estimates of the iterated kernels require a series of additional
regularization procedures absent in the classical case. Moreover,
while the classical heat kernel has a simple expression in
elementary functions, in the fractional case we have to deal with
Fox's H-functions and their asymptotics, which makes even simple
(in principle) transformations and estimates difficult technical
tasks. Note also that, in order to consider an inhomogeneous
equation (1.1), we have to develop Levi's method twice --
separately for $Z$ and $Y$.

The results can be generalized easily to the situation, in which
all the coefficients of $B$ except the leading ones depend also on
$t$. Moreover, if we consider only the case of the zero initial
functions, that is we look only for the function $Y$, then we can
include, without any changes, the general case of time-dependent
coefficients. The same is true for the case $n=1$ which resembles
the classical theory of parabolic equations. In fact, it seems
probable that the stationarity assumption can be dropped
completely. However, then the machinery of Fox's H-functions would
not be available, and we would have to deal directly with contour
integral representations of solutions of fractional diffusion
equations with coefficients depending only on $t$. Their study is
a complicated task in itself, and it looks reasonable to consider
first the stationary case, in order to identify the main
differences between the equations (1.1) and classical parabolic
equations.

The main results of this paper are collected in Sect. 2. Sect. 3
contains the information regarding H-functions, and miscellaneous lemmas
used subsequently. Proofs of parametrix estimates are given in Sect. 4
and used in Sect. 5 for substantiating the Levi method for the case
where $n\ge 2$. The case $n=1$ is considered in Sect. 6. In Sect.
7 we prove the nonnegativity of the functions $Z$ and $Y$.

\medskip
\section{MAIN RESULTS}

In this section we describe the assumptions on the coefficients
and formulate principal results. The proofs will be given in
subsequent sections.

\smallskip
{\bf 2.1.} {\it Equations with constant coefficients}. Let us
begin with the case where the coefficients of $B$ are constant,
and only the leading terms are present, so that
$$
B=\sum\limits_{i,j=1}^na_{ij}\frac{\partial^2}{\partial
x_i\partial x_j}.
$$
We assume that the matrix $A=(a_{ij})$ is positive definite.

The function $Z=Z_0(t,x-\xi )$ for this case is obtained by a
change of variables from the fundamental solution found in
\cite{SW,K2} for $B=\Delta$. In order to formulate the result, we
recall the definition of Fox's H-function which will be used
systematically in this paper.

Let $\mu ,\nu ,p,q$ be integers satisfying the conditions $0\le
\nu \le p$, $1\le \mu \le q$. Suppose we have also the complex
parameters $c_1,\ldots ,c_p$ and $d_1,\ldots ,d_q$, and positive
real parameters $\gamma_1,\ldots ,\gamma_p$ and $\delta_1,\ldots ,
\delta_q$, such that $P_1\cap P_2=\varnothing$ where
$$
P_1=\{ s=-(d_j+k)/\delta_j,\ j=1,\ldots ,\mu;\ k=0,1,2,\ldots \},
$$
$$
P_2=\{ s=(1-c_j+k)/\gamma_j,\ j=1,\ldots ,\nu;\ k=0,1,2,\ldots \}.
$$
We will need only the case where
\begin{equation}
\rho \overset{\mbox{\scriptsize def}}{=}\sum\limits_{k=1}^q\delta_k-\sum\limits_{k=1}^p
\gamma_k>0.
\end{equation}

The H-function
$$
H_{pq}^{\mu \nu}(z)=H_{pq}^{\mu \nu}\left[ z\left| \begin{matrix}
(c_1,\gamma_1), & \ldots , & (c_\nu ,\gamma_\nu ); &
(c_{\nu +1},\gamma_{\nu +1}), & \ldots , & (c_p,\gamma_p)\\
(d_1,\delta_1), & \ldots , & (d_\mu ,\delta_\mu ); &
(d_{\mu +1},\delta_{\mu +1}), & \ldots , & (d_q,\delta_q)\end{matrix}\right.
\right]
$$
is defined by the contour integral
\begin{equation}
H_{pq}^{\mu \nu}(z)=\frac{1}{2\pi
i}\int\limits_L\frac{C(s)D(s)}{E(s)F(s)}z^{-s}\,ds
\end{equation}
with
\begin{gather*}
C(s)=\prod\limits_{k=1}^\mu \Gamma (d_k+\delta_ks),\quad
D(s)=\prod\limits_{k=1}^\nu \Gamma (1-c_k-\gamma_ks),\\
E(s)=\prod\limits_{k=\mu +1}^q \Gamma (1-d_k-\delta_ks),\quad
F(s)=\prod\limits_{k=\nu +1}^p \Gamma (c_k+\gamma_ks)
\end{gather*}
($P_1$ and $P_2$ are the sets of poles for $C(s)$ and $D(s)$
respectively). The integration contour $L$ is an infinite loop
running between $s=-\infty -i\sigma$ and $s=-\infty +i\sigma$
where $\sigma >\max\limits_{1\le j\le \mu}\{ |\I d_j|/\delta_j\}$ in
such a way that $P_1$ lies to the left of $L$, and $P_2$ to the right of
$L$.

Nearly all classical special functions can be represented as H-functions
with appropriate parameters. The theory of the H-function including
its analytic properties, asymptotics, various relations, is
expounded in \cite{Br,PBM,SGG}. Below we formulate all the results
we need.

Under the above assumptions the function $H_{pq}^{\mu \nu}(z)$ is
holomorphic in a certain sector containing the positive real half-axis.

Now we can write the formula for the fundamental solution $Z_0$:
\begin{multline}
Z_0(t,x-\xi )=\frac{\pi^{-n/2}}{(\det A)^{1/2}}\left[
\sum\limits_{i,j=1}^nA^{(ij)}(x_i-\xi_i)(x_j-\xi_j)\right]^{-n/2}\\
\times H_{12}^{20}\left[ \frac{1}{4}t^{-\alpha}
\sum\limits_{i,j=1}^nA^{(ij)}(x_i-\xi_i)(x_j-\xi_j)\left|
\begin{matrix} (1,\alpha ) & \\
(\frac{n}2,1), & (1,1)\end{matrix}\right. \right]
\end{multline}
where $A=\left( A^{(ij)}\right)$ is the matrix inverse to
$(a_{ij})$. Two different proofs of (2.3) can be found (for
$B=\Delta$) in \cite{K2} and \cite{SW} (where an equivalent formula
is given). In \cite{K2} we used the Fourier transform with respect
to the spatial variables; the resulting equation is solved using
the Mittag-Leffler function, and then the inverse Fourier
transform is performed on the basis of the appropriate formulas
for H-functions. The authors of \cite{SW} used the Mellin
transform in $t$, the explicit expression of the Green function
for the Laplacian, the inverse Mellin transform, and integration
formulas for H-functions.

The function $Y=Y_0(t-\lambda ,x-y)$ appearing in the
representation (1.4), for this case has the form
\begin{multline}
Y_0(t-\lambda,x-y)=\frac{\pi^{-n/2}}{(\det A)^{1/2}}\left[
\sum\limits_{i,j=1}^nA^{(ij)}(x_i-y_i)(x_j-y_j)\right]^{-n/2}\\
\times (t-\lambda )^{\alpha -1}H_{12}^{20}\left[
\frac{1}{4}(t-\lambda)^{-\alpha}
\sum\limits_{i,j=1}^nA^{(ij)}(x_i-y_i)(x_j-y_j)\left|
\begin{matrix} (\alpha ,\alpha ) & \\
(\frac{n}2,1), & (1,1)\end{matrix}\right. \right].
\end{multline}
In fact, $Y_0(t,x)$ is the Riemann-Liouville derivative of
$Z_0(t,x)$ in $t$, of the order $1-\alpha$ (for $x\ne 0$,
$Z_0(t,x)\to 0$ as $t\to 0$, so that the Riemann-Liouville
derivative coincides in the case with the regularized fractional
derivative).

Estimates of the function $Z_0$ and $Y_0$, and of their
derivatives, are given in the following propositions. Denote
$R=t^{-\alpha }|x|^2$. Here and below the letters $C,\sigma$ will
denote various positive constants.

\medskip
\begin{prop}[see \cite{K2}]
\begin{description}
\item[(i)] If $R\ge 1$, then
\begin{equation}
\left| D_x^mZ_0(t,x)\right| \le Ct^{-\frac{\alpha (n+m)}2}\exp
\left( -\sigma t^{-\frac{\alpha }{2-\alpha
}}|x|^{\frac{2}{2-\alpha }}\right) ,\quad |m|\le 3,
\end{equation}
\begin{equation}
\left| \D Z_0(t,x)\right| \le Ct^{-\frac{\alpha (n+2)}2}\exp
\left( -\sigma t^{-\frac{\alpha }{2-\alpha
}}|x|^{\frac{2}{2-\alpha }}\right) ;
\end{equation}
\item[(ii)] If $R\le 1$, $x\ne 0$, then
\begin{equation}
\left| D_x^mZ_0(t,x)\right| \le Ct^{-\alpha}|x|^{-n+2-|m|},\quad |m|\le 3,
\end{equation}
if \ $n\ge 3$, \ or $n=2$, $m\ne 0$;
\begin{equation}
\left|Z_0(t,x)\right| \le Ct^{-\alpha}[|\log (t^{-\alpha }|x|^2)|+1],
\end{equation}
if \ $n=2$;
\begin{equation}
\left| \frac{\partial^mZ_0(t,x)}{\partial x^m}\right|
\le Ct^{-\frac{(m+1)\alpha}{2}},
\end{equation}
if \ $n=1$.
\item[(iii)] If $R\le 1$, $x\ne 0$, then
\begin{equation}
\left| \D Z_0(t,x)\right| \le \begin{cases}
Ct^{-2\alpha}|x|^{-n+2}, & \text{if $n\ge 3$},\\
Ct^{-\alpha}\left[ \left| \log (t^{-\alpha}|x|^2)\right|
+1\right], & \text{if $n=2$},\\
Ct^{-3\alpha /2}, & \text{if $n=1$}. \end{cases}
\end{equation}
\end{description}
\end{prop}

Note that the orders of the singularities at $x=0$ in (2.7),
(2.8), and (2.10) are precise (they are based on the asymptotic
expansions of the H-functions).

The next proposition contains estimates of the function $Y_0$ and
its derivatives.

\medskip
\begin{prop}
\begin{description}
\item[(i)] If $R\ge 1$, then
\begin{equation}
\left| D_x^mY_0(t,x)\right| \le Ct^{-\frac{\alpha (n+m)}2-1+\alpha}\exp
\left( -\sigma t^{-\frac{\alpha }{2-\alpha
}}|x|^{\frac{2}{2-\alpha }}\right) ,\quad |m|\le 3,
\end{equation}
\item[(ii)] If $R\le 1$, $x\ne 0$, $n>4$, then
\begin{equation}
\left| D_x^mY_0(t,x)\right| \le Ct^{-\alpha -1}|x|^{-n+4-|m|},\quad |m|\le 3,
\end{equation}
\item[(iii)] If $R\le 1$, $x\ne 0$, $n=4$, then
\begin{equation}
|Y_0(t,x)| \le Ct^{-\alpha -1}[|\log (t^{-\alpha }|x|^2)|+1],
\end{equation}
\begin{equation}
\left| D_xY_0(t,x)\right| \le Ct^{-\frac{3\alpha}{2} -1},
\end{equation}
\begin{equation}
\left| D_x^mY_0(t,x)\right| \le Ct^{-2\alpha -1}[|\log
(t^{-\alpha }|x|^2)|+1],\quad |m|=2,
\end{equation}
\begin{equation}
\left| D_x^mY_0(t,x)\right| \le Ct^{-2\alpha -1}|x|^{-1}[|\log
(t^{-\alpha }|x|^2)|+1],\quad |m|=3.
\end{equation}
\item[(iv)] If $R\le 1$, $x\ne 0$, $n=3$, then
\begin{equation}
|Y_0(t,x)| \le Ct^{-\frac{\alpha}2 -1},
\end{equation}
\begin{equation}
\left| D_xY_0(t,x)\right| \le Ct^{-\alpha -1},
\end{equation}
\begin{equation}
\left| D_x^mY_0(t,x)\right| \le Ct^{-\alpha -1}|x|^{-1},\quad
|m|=2,
\end{equation}
\begin{equation}
\left| D_x^mY_0(t,x)\right| \le Ct^{-\alpha -1}|x|^{-2},\quad
|m|=3.
\end{equation}
\item[(v)] If $R\le 1$, $x\ne 0$, $n=2$, then
\begin{equation}
|Y_0(t,x)| \le Ct^{-1},
\end{equation}
\begin{equation}
\left| D_xY_0(t,x)\right| \le Ct^{-\frac{\alpha}{2} -1},
\end{equation}
\begin{equation}
\left| D_x^mY_0(t,x)\right| \le Ct^{-\alpha -1}[|\log
(t^{-\alpha }|x|^2)|+1],\quad |m|=2,
\end{equation}
\begin{equation}
\left| D_x^mY_0(t,x)\right| \le Ct^{-\alpha -1}|x|^{-1}[|\log
(t^{-\alpha }|x|^2)|+1],\quad |m|=3.
\end{equation}
\item[(vi)] If $R\le 1$, $x\ne 0$, $n=1$, then
\begin{equation}
\left| D_x^mY_0(t,x)\right| \le Ct^{-\frac{(m-1)\alpha}2-1},\quad
m=0,1,2,3.
\end{equation}
\end{description}
\end{prop}

\medskip
{\bf 2.2.} {\it The general case}. We make the following
assumptions on the coefficients of the operator $B$.
\begin{description}
\item[(B${}_1$)] The coefficients $a_{ij}(x)$, $b_j(x)$, $c(x)$
are bounded H\"older continuous functions on $\mathbb R^n$.
\item[(B${}_2$)] The uniform parabolicity condition: there exists
such a constant $\delta >0$ that for any $x,\xi \in \mathbb R^n$
$$
\sum\limits_{i,j=1}^na_{ij}(x)\xi_i\xi_j\ge \delta |\xi |^2.
$$
\end{description}

In order to describe regularity properties of solutions of the
equation (1.1) with respect to the variable $t$, it is convenient
to introduce the class $H_\mu^{\alpha +\lambda }[0,T]$ of such
functions $\varphi (t)$ that $t^\mu \varphi (t)$ is H\"older
continuous on $[0,T]$ with the exponent $\alpha +\lambda$,
$\lambda >0$. It was proved in \cite{K2} that the problem (1.1), (1.3)
cannot have more than one bounded classical solution belonging
(for each fixed $x\in \mathbb R^n$) to $H_\mu^{\alpha +\lambda
}[0,T]$, if $0<\lambda <1-\alpha$, $0\le \mu <\lambda +1$.

Below we denote by $\gamma$ various H\"older exponents with
respect to spatial variables (without restricting generality they
will be assumed equal).

\begin{teo}
a) There exists a Green matrix $\{ Z(t,x;\xi ),Y(t,x;\xi )\}$
for the problem (1.1), (1.3), of the form
\begin{gather*}
Z(t,x;\xi )=Z_0(t,x-\xi ;\xi )+V_Z(t,x;\xi ),\\
Y(t,x;\xi )=Y_0(t,x-\xi ;\xi )+V_Y(t,x;\xi ),
\end{gather*}
where $(Z_0,Y_0)$, the Green matrix of the Cauchy problem for the
equation obtained by ``freezing'' the coefficients $a_{ij}$ at the
parameter point $\xi$ and setting other coefficients equal to
zero, satisfies the estimates listed in Propositions 1,2, with the
constants independent of $\xi \in \mathbb R^n$. The functions
$V_Z,V_Y$ satisfy the estimates
\begin{description}
\item{{\rm (i)}}
If $n=1$, then
\begin{gather*}
\left| D_x^mV_Z(t,x;\xi )\right| \le Ct^{(\gamma -m-1)\alpha
/2}\exp \left( -\sigma t^{-\frac{\alpha }{2-\alpha
}}|x-\xi |^{\frac{2}{2-\alpha }}\right) ,\\
\left| D_x^mV_Y(t,x;\xi )\right| \le Ct^{\alpha -1 +(\gamma -m-1)\alpha
/2}\exp \left( -\sigma t^{-\frac{\alpha }{2-\alpha
}}|x-\xi |^{\frac{2}{2-\alpha }}\right) ,
\end{gather*}
$m=0,1,2$.

\item{{\rm (ii)}}
If $n=2$, then
\begin{gather*}
\left| V_Z(t,x;\xi )\right| \le Ct^{\frac{\gamma \alpha }2-\alpha}
\exp \left( -\sigma t^{-\frac{\alpha }{2-\alpha
}}|x-\xi |^{\frac{2}{2-\alpha }}\right) ,\\
\left| D_xV_Z(t,x;\xi )\right| \le Ct^{\frac{\gamma_0\alpha }2-\alpha}
|x-\xi |^{-1+\gamma -\gamma_0}\exp \left(
-\sigma t^{-\frac{\alpha }{2-\alpha
}}|x-\xi |^{\frac{2}{2-\alpha }}\right) ,\\
\left| V_Y(t,x;\xi )\right| \le Ct^{\frac{\gamma \alpha }2-1}
\exp \left( -\sigma t^{-\frac{\alpha }{2-\alpha
}}|x-\xi |^{\frac{2}{2-\alpha }}\right) ,\\
\left| D_xV_Y(t,x;\xi )\right| \le Ct^{\frac{\gamma_0\alpha }2-1}
|x-\xi |^{-1+\gamma -\gamma_0}\exp \left(
-\sigma t^{-\frac{\alpha }{2-\alpha
}}|x-\xi |^{\frac{2}{2-\alpha }}\right) ,
\end{gather*}
where $\gamma_0<\gamma$ is an arbitrary fixed positive constant
(here and below in the formulas containing $\gamma_0$ the constant
$C$ depends on $\gamma_0$).

\item{{\rm (iii)}}
If $n=3$, then
\begin{gather*}
\left| D_x^mV_Z(t,x;\xi )\right| \le Ct^{\frac{\gamma_0\alpha }2-\alpha}
|x-\xi |^{-1-|m|+\gamma -\gamma_0}\exp
\left( -\sigma t^{-\frac{\alpha }{2-\alpha
}}|x-\xi |^{\frac{2}{2-\alpha }}\right) ,\\
\left| D_x^mV_Y(t,x;\xi )\right| \le Ct^{(\gamma_0+\gamma )\frac{\alpha}4-1}
|x-\xi |^{-1-|m|+(\gamma -\gamma_0)/2}
\exp \left( -\sigma t^{-\frac{\alpha }{2-\alpha
}}|x-\xi |^{\frac{2}{2-\alpha }}\right) ,
\end{gather*}
$|m|=0,1$.

\item{{\rm (iv)}}
If $n=4$, then
\begin{gather*}
\left| D_x^mV_Z(t,x;\xi )\right|
\le Ct^{\frac{(\gamma -\gamma_0)\alpha }2-\alpha}
|x-\xi |^{-2-|m|+\gamma -2\gamma_0}\exp \left( -\sigma
t^{-\frac{\alpha }{2-\alpha }}|x-\xi |^{\frac{2}{2-\alpha }}\right) ,\\
\left| D_x^mV_Y(t,x;\xi )\right| \le Ct^{\frac{(\gamma -\gamma_0)\alpha }4-1}
|x-\xi |^{-2-|m|+\gamma -2\gamma_0}
\exp \left( -\sigma t^{-\frac{\alpha }{2-\alpha
}}|x-\xi |^{\frac{2}{2-\alpha }}\right) ,
\end{gather*}
$|m|=0,1$; $\gamma_0<\gamma /2$.

\item{{\rm (v)}}
If $n\ge 5$, then
\begin{gather*}
\left| D_x^mV_Z(t,x;\xi )\right| \le Ct^{\frac{\gamma_0\alpha }2-\alpha}
|x-\xi |^{2-n-|m|+\gamma -\gamma_0}\exp
\left( -\sigma t^{-\frac{\alpha }{2-\alpha
}}|x-\xi |^{\frac{2}{2-\alpha }}\right) ,\\
\left| D_x^mV_Y(t,x;\xi )\right| \le Ct^{(\gamma_0+\gamma )\frac{\alpha}4-1}
|x-\xi |^{2-n-|m|+(\gamma -\gamma_0)/2}
\exp \left( -\sigma t^{-\frac{\alpha }{2-\alpha
}}|x-\xi |^{\frac{2}{2-\alpha }}\right) ,
\end{gather*}
$|m|=0,1$.

\item{{\rm (vi)}}
If $n\ge 2,|m|=2$, then
\begin{gather*}
\left| D_x^mV_Z(t,x;\xi )\right| \le Ct^{\gamma_1-\alpha }
|x-\xi |^{-n+\gamma_2}\exp \left( -\sigma t^{-\frac{\alpha }{2-\alpha
}}|x-\xi |^{\frac{2}{2-\alpha }}\right) ,\\
\left| D_x^mV_Y(t,x;\xi )\right| \le Ct^{\gamma_1-1}|x-\xi |^{-n+\gamma_2}
\exp \left( -\sigma t^{-\frac{\alpha }{2-\alpha
}}|x-\xi |^{\frac{2}{2-\alpha }}\right) ,
\end{gather*}
where $\gamma_1,\gamma_2$ are some positive constants.

\item{{\rm (vii)}}
The estimates for the fractional derivatives $\D V_Z$, $\D V_Y$ of
the functions $Z,Y$ are the same as the above estimates for their
second order derivatives.
\end{description}

b) Let $u(t,x)$ be the solution (1.4) of the problem (1.1), (1.3)
with a bounded continuous function $u_0$ (locally H\"older
continuous if $n>1$) and a bounded jointly continuous and
locally H\"older continuous, in $x$,
function $f$. If $0<\lambda <1-\alpha$, $\alpha +\lambda <\nu$,
then for any $x\in \mathbb R^n$ $u(\cdot ,x)\in H_\nu ^{\alpha
+\lambda}[0,T]$.

c) The functions $Z$ and $Y$ are nonnegative.
\end{teo}

\medskip
The rest of this paper is devoted to the proof of the Theorem. As
it was mentioned in the Introduction, it can be generalized to the
nonstationary case, and some lemmas below are formulated in a way
general enough to cover that more general situation.

\medskip
\section{AUXILIARY RESULTS}

{\bf 3.1.} {\it H-functions}. We collect here some results
regarding H-functions, which will be used below. For the proofs
and further details see \cite{Br,PBM,SGG}.

\smallskip
a) {\it Differentiation formulas}.
\begin{equation}
\frac{d}{dz}H_{12}^{20}\left[ z\left| \begin{matrix}
(c_1,\gamma_1) & \\
(d_1,\delta_1), & (d_2,\delta_2)\end{matrix} \right. \right]
=-z^{-1}H_{23}^{30}\left[ z\left| \begin{matrix}
(c_1,\gamma_1), & (0,1) & \\
(d_1,\delta_1), & (d_2,\delta_2), & (1,1)\end{matrix}\right. \right] .
\end{equation}
\begin{equation}
\frac{d}{dz}H_{23}^{30}\left[ z\left| \begin{matrix}
(c_1,\gamma_1), & (c_2,\gamma_2) & \\
(d_1,\delta_1), & (d_2,\delta_2), & (d_3,\delta_3)\end{matrix}
\right. \right]
=-z^{-1}H_{34}^{40}\left[ z\left| \begin{matrix}
(c_1,\gamma_1), & (c_2,\gamma_2), & (0,1) & \\
(d_1,\delta_1), & (d_2,\delta_2), & (d_3,\gamma_3) & (1,1)\end{matrix}
\right. \right] .
\end{equation}

In a similar way $\frac{d}{dz}H_{34}^{40}$ is expressed via
$H_{45}^{50}$ etc.

By (2.2),
\begin{equation}
H_{12}^{20}\left[ z\left| \begin{matrix}
(1,\alpha ) & \\
(\frac{n}2,1), & (1,1)\end{matrix}\right. \right] =\frac{1}{2\pi
i}\int\limits_L\frac{\Gamma (\frac{n}2+s)\Gamma (1+s)}{\Gamma
(1+\alpha s)}z^{-s}\,ds.
\end{equation}
We will need this function with $z=\omega (t-\tau )^{-\alpha }$,
$\omega >0$, $t>\tau$. We have
\begin{equation}
u(t)=H_{12}^{20}\left[ \omega (t-\tau )^{-\alpha }\left| \begin{matrix}
(1,\alpha ) & \\
(\frac{n}2,1), & (1,1)\end{matrix}\right. \right] =\frac{1}{2\pi
i}\int\limits_L\frac{\Gamma (\frac{n}2+s)\Gamma (1+s)}{\Gamma
(1+\alpha s)}\omega^{-s}(t-\tau )^{\alpha s}\,ds.
\end{equation}
In order to calculate the Riemann-Liouville fractional derivative
$\left( D_{\tau +}^\alpha u\right) (t)=\frac{1}{\Gamma (1-\alpha
)}\frac{\partial }{\partial t}\int\limits_\tau ^t(t-s)^{-\alpha
}u(s)\,ds$, we note that $D_{\tau +}^\alpha$ transforms $(t-\tau )^{\alpha
s}$ into $\dfrac{\Gamma (1+\alpha s)}{\Gamma (1-\alpha +\alpha
s)}(t-\tau )^{\alpha s-\alpha }$ (see \cite{SKM}). Now we find
from (3.4) that
\begin{equation}
\left( D_{\tau +}^\alpha u\right) (t)=(t-\tau )^{-\alpha }
H_{12}^{20}\left[ \omega (t-\tau )^{-\alpha }\left| \begin{matrix}
(1-\alpha ,\alpha ) & \\
(\frac{n}2,1), & (1,1)\end{matrix}\right. \right]
\end{equation}

Similarly
\begin{equation}
\left( D_{\tau +}^{1-\alpha }u\right) (t)=(t-\tau )^{\alpha -1}
H_{12}^{20}\left[ \omega (t-\tau )^{-\alpha }\left| \begin{matrix}
(\alpha ,\alpha ) & \\
(\frac{n}2,1), & (1,1)\end{matrix}\right. \right]
\end{equation}

\medskip
b) {\it Asymptotics at infinity}. The asymptotic behavior of
H-functions for $z\to \infty$ has been thoroughly investigated. We
will need only one result of this kind, for a specific class of
H-functions, for a real argument, with only the leading term of
the asymptotic expansion:
\begin{equation}
H_{p\mu}^{\mu 0}\left[ z\left| \begin{matrix}
(c_1,\gamma_1), & \ldots , &  (c_p,\gamma_p)\\
(d_1,\delta_1), & \ldots , & (d_\mu ,\delta_\mu )\end{matrix}\right.
\right] \sim \text{const}\cdot z^{(1-a)/\rho}\exp \left(
-b^{1/\rho}\rho z^{1/\rho}\right)
\end{equation}
where $a=\sum\limits_{k=1}^pc_k-\sum\limits_{k=1}^\mu
d_k+\frac{1}2(\mu -p+1)$,
$b=\prod\limits_{k=1}^p\gamma_k^{\gamma_k}\prod\limits_{l=1}^\mu
\delta_l^{-\delta_l}$, and $\rho$ is given by (2.1).

In particular, from (3.7) we get estimates of some specific
H-functions for $|z|\ge 1$:
\begin{equation}
\left| H_{12}^{20}\left[ z\left| \begin{matrix}
(1,\alpha ) & \\
(\frac{n}2,1), & (1,1)\end{matrix}\right. \right]\right| \le
C|z|^{\frac{n}{2(2-\alpha )}}\exp \left( -\sigma
|z|^{\frac{1}{2-\alpha}}\right);
\end{equation}
\begin{equation}
\left| H_{12}^{20}\left[ z\left| \begin{matrix}
(\alpha ,\alpha ) & \\
(\frac{n}2,1), & (1,1)\end{matrix}\right. \right]\right| \le
C|z|^{\frac{n+2-2\alpha }{2(2-\alpha )}}\exp \left( -\sigma
|z|^{\frac{1}{2-\alpha}}\right).
\end{equation}
Similar estimates hold for functions appearing in expressions for
derivatives of the above H-functions. In all cases the estimates
contain the same exponentially decreasing factor while the degrees
of positive powers of $|z|$ are different and specific for each
case.

\medskip
c) {\it Asymptotics near the origin}. In order to obtain, for the
above class of H-functions, an asymptotic expansion near the
origin, we write an H-function as
$$
H_{p\mu}^{\mu 0}\left[ z\left| \begin{matrix}
(c_1,\gamma_1), & \ldots , &  (c_p,\gamma_p)\\
(d_1,\delta_1), & \ldots , & (d_\mu ,\delta_\mu )\end{matrix}\right.
\right] =\sum\limits_{s\in P_1}\Res \left(
\frac{C(s)}{F(s)}z^{-s}\right).
$$

In specific cases, writing $P_1$ explicitly and using well-known
properties of the function Gamma, we find a required number of
terms in asymptotic expansions. In particular, for $z\to +0$
\begin{equation}
H_{12}^{20}\left[ z\left| \begin{matrix}
(1,\alpha ) & \\
(\frac{n}2,1), & (1,1)\end{matrix}\right. \right] =\begin{cases}
\frac{\Gamma (\frac{n}2-1)}{\Gamma (1-\alpha )}z+O\left( z^{\min
(\frac{n}2,2)}\right) , & \text{if $n\ge 3$},\\
-\frac{z\log z}{\Gamma (1-\alpha )}+O(z), & \text{if $n=2$},\\
\frac{\sqrt{\pi }}{\Gamma (1-\frac{\alpha}2)}z^{1/2}-
\frac{2\sqrt{\pi }}{\Gamma (1-\alpha)}z+\frac{2\sqrt{\pi
}}{\Gamma (1-\frac{3\alpha}2)}z^{3/2}+O(z^2), & \text{if $n=1$}
\end{cases}
\end{equation}
(here and below we present those terms of the asymptotic
expansions which are actually used in this work and in \cite{K2}),
\begin{equation}
H_{23}^{30}\left[ z\left| \begin{matrix}
(1,\alpha ), & (0,1), &  \\
(\frac{n}2,1), & (1,1), & (1,1)\end{matrix}\right. \right] =\begin{cases}
-\frac{\Gamma (\frac{n}2-1)}{\Gamma (1-\alpha )}z+O\left( z^{\min
(\frac{n}2,2)}\right) , & \text{if $n\ge 3$},\\
\frac{z\log z}{\Gamma (1-\alpha )}+O(z), & \text{if $n=2$},\\
-\frac{\sqrt{\pi }}{2\Gamma (1-\frac{\alpha}2)}z^{1/2}+
\frac{2\sqrt{\pi }}{\Gamma (1-\alpha)}z-\frac{3\sqrt{\pi
}}{\Gamma (1-\frac{3\alpha}2)}z^{3/2}+O(z^2), & \text{if $n=1$},
\end{cases}
\end{equation}
\begin{equation}
H_{34}^{40}\left[ z\left| \begin{matrix}
(1,\alpha ), & (0,1), & (0,1) & \\
(\frac{n}2,1), & (1,1), & (1,1), & (1,1)\end{matrix}\right.
\right] =\begin{cases}
\frac{\Gamma (\frac{n}2-1)}{\Gamma (1-\alpha )}z+O\left( z^{\min
(\frac{n}2,2)}\right) , & \text{if $n\ge 3$},\\
-\frac{z\log z}{\Gamma (1-\alpha )}+O(z), & \text{if $n=2$},\\
\frac{\sqrt{\pi }}{4\Gamma (1-\frac{\alpha}2)}z^{1/2}-
\frac{2\sqrt{\pi }}{\Gamma (1-\alpha)}z\\ +\frac{9\sqrt{\pi
}}{2\Gamma (1-\frac{3\alpha}2)}z^{3/2}+O(z^2), & \text{if $n=1$},
\end{cases}
\end{equation}
\begin{equation}
H_{45}^{50}\left[ z\left| \begin{matrix}
(1,\alpha ), & (0,1), & (0,1), & (0,1) & \\
(\frac{n}2,1), & (1,1), & (1,1), & (1,1), & (1,1)\end{matrix}\right.
\right] =\begin{cases}
-\frac{\Gamma (\frac{n}2-1)}{\Gamma (1-\alpha )}z+O\left( z^{\min
(\frac{n}2,2)}\right) , & \text{if $n\ge 3$},\\
\frac{z\log z}{\Gamma (1-\alpha )}+O(z), & \text{if $n=2$},\\
-\frac{\sqrt{\pi }}{8\Gamma (1-\frac{\alpha}2)}z^{1/2}+
\frac{2\sqrt{\pi }}{\Gamma (1-\alpha)}z\\ -\frac{27\sqrt{\pi
}}{4\Gamma (1-\frac{3\alpha}2)}z^{3/2}+O(z^2), & \text{if $n=1$},
\end{cases}
\end{equation}

Next we give the expansions for the function appearing in (2.4)
and those emerging in the course of differentiating (2.4). A
different character of their asymptotics is caused by the fact
that in this case the pole $s=1$ appears both in the numerator and
denominator of the fraction in (2.2) and cancel each other. We
have
\begin{equation}
H_{12}^{20}\left[ z\left| \begin{matrix}
(\alpha ,\alpha ) & \\
(\frac{n}2,1), & (1,1)\end{matrix}\right. \right] =\begin{cases}
-\frac{\Gamma (\frac{n}2-2)}{\Gamma (-\alpha )}z^2+O\left( z^{\min
(\frac{n}2,3)}\right) , & \text{if $n\ge 5$},\\
\frac{z^2\log z}{\Gamma (-\alpha )}+\gamma_1z^2+O(z^3\log z),
& \text{if $n=4$},\\
-\frac{2\sqrt{\pi }}{\Gamma (-\frac{\alpha}2)}z^{3/2}+
\frac{2\sqrt{\pi }}{\Gamma (-\alpha)}z^2+O(z^{5/2}), & \text{if
$n=3$},\\
\alpha z-\frac{z^2\log z}{\Gamma (-\alpha )}+\gamma_2z^2+O(z^3\log z),
& \text{if $n=2$},\\
\frac{\sqrt{\pi }}{\Gamma (\frac{\alpha}2)}z^{1/2}+
\frac{2\sqrt{\pi }}{\Gamma (-\frac{\alpha}2)}z^{3/2}+O(z^2), & \text{if
$n=1$},\end{cases}
\end{equation}
where $\gamma_1,\gamma_2$ are certain constants. Next,
\begin{equation}
H_{23}^{30}\left[ z\left| \begin{matrix}
(\alpha ,\alpha ), & (0,1) & \\
(\frac{n}2,1), & (1,1), & (1,1)\end{matrix}\right. \right] =\begin{cases}
\frac{2\Gamma (\frac{n}2-2)}{\Gamma (-\alpha )}z^2+O\left( z^{\min
(\frac{n}2,3)}\right) , & \text{if $n\ge 5$},\\
-2\frac{z^2\log z}{\Gamma (-\alpha )}-2\gamma_1z^2+O(z^3\log z),
& \text{if $n=4$},\\
\frac{3\sqrt{\pi }}{\Gamma (-\frac{\alpha}2)}z^{3/2}-
\frac{4\sqrt{\pi }}{\Gamma (-\alpha)}z^2+O(z^{5/2}), & \text{if
$n=3$},\\
-\alpha z-2\frac{z^2\log z}{\Gamma (-\alpha )}-2\gamma_2z^2+O(z^3\log z),
& \text{if $n=2$},\\
-\frac{\sqrt{\pi }}{2\Gamma (\frac{\alpha}2)}z^{1/2}-
\frac{3\sqrt{\pi }}{\Gamma (-\frac{\alpha}2)}z^{3/2}+O(z^2), & \text{if
$n=1$},\end{cases}
\end{equation}
\begin{equation}
H_{34}^{40}\left[ z\left| \begin{matrix}
(\alpha ,\alpha ), & (0,1), & (0,1) &\\
(\frac{n}2,1), & (1,1), & (1,1), & (1,1)\end{matrix}\right.
\right] =\begin{cases}
-\frac{4\Gamma (\frac{n}2-2)}{\Gamma (-\alpha )}z^2+O\left( z^{\min
(\frac{n}2,3)}\right) , & \text{if $n\ge 5$},\\
4\frac{z^2\log z}{\Gamma (-\alpha )}+4\gamma_1z^2+O(z^3\log z),
& \text{if $n=4$},\\
-\frac{9\sqrt{\pi }}{2\Gamma (-\frac{\alpha}2)}z^{3/2}+
\frac{8\sqrt{\pi }}{\Gamma (-\alpha)}z^2+O(z^{5/2}), & \text{if
$n=3$},\\
\alpha z-4\frac{z^2\log z}{\Gamma (-\alpha )}+4\gamma_2z^2+O(z^3\log z),
& \text{if $n=2$},\\
\frac{\sqrt{\pi }}{4\Gamma (\frac{\alpha}2)}z^{1/2}+
\frac{9\sqrt{\pi }}{2\Gamma (-\frac{\alpha}2)}z^{3/2}+O(z^2), & \text{if
$n=1$},\end{cases}
\end{equation}
Finally,
\begin{multline}
H_{45}^{50}\left[ z\left| \begin{matrix}
(\alpha ,\alpha ), & (0,1), & (0,1), & (0,1) & \\
(\frac{n}2,1), & (1,1), & (1,1), & (1,1), & (1,1)\end{matrix}\right.
\right] \\
=\begin{cases}
\frac{8\Gamma (\frac{n}2-2)}{\Gamma (-\alpha )}z^2+O\left( z^{\min
(\frac{n}2,3)}\right) , & \text{if $n\ge 5$},\\
-8\frac{z^2\log z}{\Gamma (-\alpha )}-8\gamma_1z^2+O(z^3\log z),
& \text{if $n=4$},\\
\frac{27\sqrt{\pi }}{4\Gamma (-\frac{\alpha}2)}z^{3/2}-
\frac{16\sqrt{\pi }}{\Gamma (-\alpha)}z^2+O(z^{5/2}), & \text{if
$n=3$},\\
-\alpha z+8\frac{z^2\log z}{\Gamma (-\alpha )}-8\gamma_2z^2
+O(z^3\log z), & \text{if $n=2$},\\
-\frac{\sqrt{\pi }}{8\Gamma (\frac{\alpha}2)}z^{1/2}-
\frac{27\sqrt{\pi }}{4\Gamma (-\frac{\alpha}2)}z^{3/2}+O(z^2), & \text{if
$n=1$},\end{cases}
\end{multline}

\medskip
{\bf 3.2.} {\it Miscellaneous lemmas}. Let us fix $\beta \in
(0,1)$. For any $x,y\in \mathbb R^n$, $\lambda <t$ denote
$$
\rho (t,x;\lambda ,y)=\left( \frac{|x-y|}{(t-\lambda
)^\beta}\right)^{\frac{1}{1-\beta}}.
$$
\begin{lem}
Suppose that $0\le \tau <\lambda <t\le T$, $x,y,\xi \in \mathbb
R^n$. Then
\begin{equation}
\rho (t,x;\lambda ,y)+\rho (\lambda ,y;\tau ,\xi )\ge \rho
(t,x;\tau ,\xi ).
\end{equation}
\end{lem}

{\it Proof}. Let $f(y)=\rho (t,x;\lambda ,y)+\rho (\lambda ,y;\tau ,\xi
)$. Introducing a new variable $\widehat{y}=y-x$ we can write
$$
f(y)=f_1(\widehat{y})=\left( \frac{|\widehat{y}|}{(t-\lambda
)^\beta}\right)^{\frac{1}{1-\beta}}+\left( \frac{|\widehat{y}-(\xi
-x)|}{(\lambda -\tau )^\beta}\right)^{\frac{1}{1-\beta}}.
$$
Denoting
$$
f_2(r)=\left( \frac{r}{(t-\lambda
)^\beta}\right)^{\frac{1}{1-\beta}}+\left( \frac{|r-|\xi
-x||}{(\lambda -\tau )^\beta}\right)^{\frac{1}{1-\beta}},\quad
r\ge 0,
$$
we find that
\begin{equation}
\min\limits_{y\in \mathbb R^n}f(y)=\min\limits_{y\in \mathbb
R^n}f_1(y)\ge \min\limits_{r\ge 0}f_2(r).
\end{equation}

The derivative $f'_2(r)$ is positive for $r>|\xi -x|$; thus
$f_2(r)$ attains its minimum value on the interval $[0,|x-\xi |]$.
A standard search for its possible local extremum points yields a
single point, in which the function takes the value $\left(
\dfrac{|\xi -x|}{(t-\tau )^\beta}\right)^{\frac{1}{1-\beta}}$. On
the other hand,
$$
f_2(0)=\left( \frac{|\xi -x|}{(\lambda -\tau
)^\beta}\right)^{\frac{1}{1-\beta}},\quad f_2(|x-\xi |)=\left(
\frac{|\xi -x|}{(t-\lambda )^\beta}\right)^{\frac{1}{1-\beta}},
$$
whence $f_2(r)\ge \left(
\dfrac{|\xi -x|}{(t-\tau )^\beta}\right)^{\frac{1}{1-\beta}}$.
Together with (3.19) this implies (3.18). $\quad \Box$

\medskip
The next three lemmas establish integral inequalities involving
the function $\rho$. We keep the notation of Lemma 1.

\begin{lem}
For any $\mu >0, \varepsilon \in (0,1)$
\begin{multline}
[(t-\lambda )(\lambda -\tau )]^{-\beta n}\int\limits_{\mathbb
R^n}\exp \{-\mu [\rho (t,x;\lambda ,y)+\rho (\lambda ,y;\tau ,\xi
)]\} \\
\le C(\varepsilon )(t-\tau )^{-n\beta }\exp \{-(1-\varepsilon
)\mu \rho (t,x;\tau ,\xi )\}.
\end{multline}
\end{lem}

This lemma is a consequence of Lemma 1. See the proof of Lemma 5.1
in \cite{E} where a similar result was obtained for $n=1$.

\begin{lem}
For any $a_1,a_2>0$, any $b_1,b_2>0$, such that $b_1+b_2<n$, and
any $\mu >0$, $\varepsilon \in (0,1)$
\begin{multline}
\int\limits_\tau^t(t-\lambda )^{a_1-1}(\lambda -\tau
)^{a_2-1}\int\limits_{\mathbb R^n}|x-y|^{b_1-n}|y-\xi
|^{b_2-n}\exp \{-\mu [\rho (t,x;\lambda ,y)+\rho (\lambda ,y;\tau ,\xi
)]\}\,dy\,d\lambda \\
\le C(\varepsilon,a_1,a_2,b_1,b_2)(t-\tau )^{a_1+a_2-1}|x-\xi
|^{b_1+b_2-n}\exp \{-(1-\varepsilon )\mu \rho (t,x;\tau ,\xi )\}.
\end{multline}
\end{lem}

{\it Proof}. Writing $\mathbb R^n=\bigcup\limits_{k=1}^4V_k$ where
\begin{gather*}
V_1=\left\{ y\in \mathbb R^n:\ |x-y|>(t-\lambda )^\beta ,
\ |y-\xi |>(\lambda -\tau )^\beta \right\},\\
V_2=\left\{ y\in \mathbb R^n:\ |x-y|>(t-\lambda )^\beta ,
\ |y-\xi |\le (\lambda -\tau )^\beta \right\},\\
V_3=\left\{ y\in \mathbb R^n:\ |x-y|\le (t-\lambda )^\beta ,
\ |y-\xi |>(\lambda -\tau )^\beta \right\},\\
V_1=\left\{ y\in \mathbb R^n:\ |x-y|\le (t-\lambda )^\beta ,
\ |y-\xi |\le (\lambda -\tau )^\beta \right\},
\end{gather*}
we decompose the left-hand side in (3.21) into the sum of four
integrals $J_1,J_2,J_3,J_4$.

If $y\in V_1$, then $|x-y|^{b_1-n}<(t-\lambda )^{\beta (b_1-n)}$,
$|y-\xi |^{b_2-n}<(\lambda -\tau)^{\beta (b_2-n)}$, so that
\begin{multline*}
J_1\le \int\limits_\tau^t(t-\lambda )^{a_1+\beta b_1-1}(\lambda -\tau
)^{a_2+\beta b_2-1}[(t-\lambda )(\lambda -\tau )]^{-n\beta }\\
\times
\int\limits_{\mathbb R^n}\exp \{-\mu [\rho (t,x;\lambda ,y)+\rho
(\lambda ,y;\tau ,\xi )]\}\,dy\,d\lambda ,
\end{multline*}
and by Lemma 2
\begin{equation}
J_1\le C(t-\tau )^{a_1+a_2+(b_1+b_2-n)\beta -1}\exp
\{-(1-\frac{\varepsilon }2)\mu \rho (t,x;\tau ,\xi )\}.
\end{equation}
Since
\begin{multline*}
(t-\tau )^{(b_1+b_2-n)\beta }\exp \{-\frac{\varepsilon }2\mu \rho
(t,x;\tau ,\xi )\}\\
=[\rho (t,x;\tau ,\xi )]^{(1-\beta )(n-b_1-b_2)}\exp
\{-\frac{\varepsilon }2\mu \rho (t,x;\tau ,\xi )\}|x-\xi
|^{b_1+b_2-n}\le C|x-\xi |^{b_1+b_2-n},
\end{multline*}
the inequality (3.22) means that $J_1$ does not exceed the
right-hand side of (3.21).

In order to estimate $J_2$, we consider two distinct cases.
Suppose first that $|x-\xi |\ge 2(t-\tau )^\beta$. Then $|x-y|\ge (t-\tau
)^\beta$, so that by Lemma 1
\begin{multline}
J_2\le (t-\tau )^{\beta (b_1-n)}\exp \{-\mu \rho (t,x;\tau ,\xi )\}
\int\limits_\tau^t(t-\lambda )^{a_1-1}(\lambda -\tau
)^{a_2-1}\,d\lambda \int\limits_{|y-\xi |\le (\lambda -\tau
)^\beta }|y-\xi |^{b_2-n}\,dy\\
=C(t-\tau )^{a_1+a_2-1}(t-\tau )^{\beta (b_1+b_2-n)}\exp \{-\mu \rho
(t,x;\tau ,\xi )\}\\
\le C(\varepsilon )(t-\tau )^{a_1+a_2-1}|x-\xi
|^{b_1+b_2-n}\exp \{-(1-\varepsilon )\mu \rho (t,x;\tau ,\xi )\}.
\end{multline}

If $|x-\xi |<2(t-\tau )^\beta$, then $|x-y|<3(t-\tau )^\beta$,
$$
J_2\le \exp \{-\mu \rho (t,x;\tau ,\xi )\}\int\limits_\tau^t(t-\lambda
)^{a_1-1}(\lambda -\tau )^{a_2-1}\, d\lambda
\int\limits_{|y-\xi |<3(t-\tau )^\beta }|x-y|^{b_1-n}|y-\xi
|^{b_2-n}\, dy,
$$
and by Lemma 2 from Chapter 1 of \cite{F} we get
$$
J_2\le C(t-\tau )^{a_1+a_2-1}|x-\xi |^{b_1+b_2-n}\exp \{-\mu \rho
(t,x;\tau ,\xi )\}.
$$
Together with (3.23) this yields the required bound for $J_2$.

The estimate for $J_3$ is proved similarly to that of $J_2$.

Let us consider $J_4$. It follows from the definition of $V_4$
that $J_4=0$ if $|x-\xi |\ge 2(t-\tau )^\beta$. Therefore we may
assume that $|x-\xi |<2(t-\tau )^\beta$ and note that
$$
V_4\subset \left\{ y\in \mathbb R^n:\ |y-\xi |<2(t-\tau )^\beta
\right\}.
$$
Then
$$
J_4\le \exp \{-\mu \rho (t,x;\tau ,\xi )\}\int\limits_\tau^t(t-\lambda
)^{a_1-1}(\lambda -\tau )^{a_2-1}\, d\lambda \int\limits_{|y-\xi |<2(t-\tau
)^\beta }|x-y|^{b_1-n}|y-\xi |^{b_2-n}\, dy.
$$
Using again Lemma 2 from Chapter 1 of \cite{F} we find that
$$
J_4\le C(t-\tau )^{a_1+a_2-1}|x-\xi |^{b_1+b_2-n}\exp \{-\mu \rho
(t,x;\tau ,\xi )\},
$$
as desired. $\quad \Box$

\medskip
{\it Remark 1}. If $b_1+b_2>n$, then the assertion of Lemma 3
still holds if the factor $|x-\xi |^{b_1+b_2-n}$ is omitted from
the right-hand side of (3.21). The proof is similar; we have only
to use the appropriate estimate from the same lemma of \cite{F}.

\begin{lem}
If $0\le c_2<c_1<1$, $\gamma >0$, then
\begin{multline}
\int\limits_{\mathbb R^n}|x-y|^{-n+\gamma }\exp \{ -\mu [c_1\rho
(t,x;\lambda ,y)+c_2\rho (\lambda ,y;\tau ,\xi )]\}\,dy\\
\le C \exp \{-\mu c_2\rho (t,x;\tau ,\xi )\}(t-\lambda )^{\beta
\gamma }
\end{multline}
where the constant depends on $c_2-c_1$.
\end{lem}

{\it Proof}. Let us write the integral in the left-hand side of
(3.24) as the sum $I_1+I_2$ of two integrals corresponding to the
decomposition $\mathbb R^n=W_1\cup W_2$ where $W_1=\{y\in \mathbb
R^n:\ |y-x|\le (t-\lambda )^\beta \}$, $W_2=\{y\in \mathbb
R^n:\ |y-x|>(t-\lambda )^\beta \}$.

Using Lemma 1, we get
$$
I_1\le \exp \{-\mu c_2\rho (t,x;\tau ,\xi )\}\int\limits_{|y-x|\le
(t-\lambda )^\beta }|x-y|^{-n+\gamma }\,dy
=C(t-\lambda )^{\beta \gamma }\exp \{-\mu c_2\rho (t,x;\tau ,\xi
)\}.
$$

Next, again by Lemma 1,
\begin{multline*}
I_2=\int\limits_{|y-x|>(t-\lambda )^\beta}|x-y|^{-n+\gamma }\exp
\{ -\mu [c_1\rho
(t,x;\lambda ,y)+c_2\rho (\lambda ,y;\tau ,\xi )]\}\,dy\\
\le (t-\lambda )^{-n\beta +\gamma \beta }\exp \{-\mu c_2\rho (t,x;\tau ,
\xi )\}\int\limits_{\mathbb R^n}\exp \{-\mu (c_1-c_2)\rho (t,x;
\lambda ,y)\,dy\}\\
\le (t-\lambda )^{\beta \gamma }\exp \{-\mu c_2\rho (t,x;\tau ,\xi )\},
\end{multline*}
as desired. $\quad \Box$

\medskip
The next two lemmas give estimates of the iterated kernels
\begin{equation}
K_m(t,x;\tau ,\xi )=\int\limits_\tau^td\lambda \int\limits_{\mathbb R^n}
K_1(t,x;\lambda ,y)K_{m-1}(\lambda ,y;\tau ,\xi )\,dy,\quad m\ge 2,
\end{equation}
where $K_1=K$ is a given kernel. We will treat the cases $n\ge 2$ and
$n=1$ separately.

\begin{lem}
Let $n\ge 2$. Suppose that $K(t,x;\tau ,\xi )$ is a continuous
function on
$$
\mathcal P_n=\left\{ (t,x;\tau ,\xi ):\ x,\xi \in \mathbb R^n,x\ne \xi,
0\le \tau <t\le T\right\},
$$
such that
\begin{equation}
|K(t,x;\tau ,\xi )|\le A_1(t-\tau )^{\nu_0\beta -1}|x-\xi
|^{-n+\nu_1}\exp \{-\mu \rho (t,x;\tau ,\xi )\}
\end{equation}
where $A_1,\mu>0$, $0<\nu_0,\nu_1,\beta <1$, and $n\nu_1^{-1}$
is not an integer.
Then the series $\mathcal R(t,x;\tau ,\xi )=\sum\limits_{m=1}^\infty
K_m(t,x;\tau ,\xi
)$ is absolutely and uniformly convergent on $\mathcal P_n\cap
\{(t,x;\tau ,\xi ):\ t-\tau \ge \delta >0,|x-\xi |\ge \delta
>0\}$, and
\begin{equation}
|\mathcal R(t,x;\tau ,\xi )|\le C(t-\tau )^{\nu_0\beta -1}|x-\xi
|^{-n+\nu_1}\exp \{-\mu^*\rho (t,x;\tau ,\xi )\},
\end{equation}
with some $\mu^*\in (0,\mu )$.
\end{lem}

{\it Proof}. The bounds for the kernels (3.25) are obtained by two
stages. First we use Lemma 3. We find, for a small $\varepsilon
>0$, that
$$
|K_2(t,x;\tau ,\xi )|\le C_2(\varepsilon )(t-\tau )^{2\nu_0\beta
-1}|x-\xi |^{-n+2\nu_1}\exp \{-\mu (1-\varepsilon )\rho (t,x;\tau ,\xi )\},
$$
and so on, so that
\begin{equation}
|K_s(t,x;\tau ,\xi )|\le C_s(\varepsilon )(t-\tau )^{s\nu_0\beta
-1}|x-\xi |^{-n+s\nu_1}\exp \{-\mu (1-(s-1)\varepsilon )\rho (t,x;\tau ,\xi
)\},
\end{equation}
whenever $-n+s\nu_1<0$. Let $s^*=\min \{s:\ -n+s\nu_1>0\}$. Then
by Remark 1
\begin{equation}
|K_{s^*}(t,x;\tau ,\xi )|\le C_0(t-\tau )^{a_0-1}\exp \{-\mu^*
\rho (t,x;\tau ,\xi )\},
\end{equation}
where $a_0=s^*\nu_0\beta$, $\mu^*=[1-(s^*-1)\varepsilon ]\mu$.

For all the next iterations, we use (3.26), (3.29), and Lemma 4.
Using the identity
$$
\int\limits_\tau^t(t-\lambda )^{a_1-1}(\lambda -\tau
)^{a_2-1}\,d\lambda =B(a_1,a_2)(t-\tau )^{a_1+a_2-1}
$$
we get by induction that
\begin{multline}
|K_{s^*+m}(t,x;\tau ,\xi )|\le C_0C_1^m\left\{
\prod_{j=1}^mB((\nu_0+\nu_1)\beta
,a_0+(\nu_0+\nu_1)(j-1))\right\}\\ \times
(t-\tau )^{m(\nu_0+\nu_1)\beta +a_0}\exp \{-\mu^*\rho (t,x;\tau ,\xi
)\}.
\end{multline}
The product in the right-hand side of (3.30) equals $\dfrac{\Gamma
(a_0)[\Gamma ((\nu_0+\nu_1)\beta )]^{m-1}}{\Gamma
(a_0+(m-1)(\nu_0+\nu_1)\beta )}$. This implies the convergence of
the series for $\mathcal R$. The estimate (3.27) follows from (3.28) and
(3.30). $\quad \Box$

\medskip
Let us consider the iterations (3.25) for $n=1$, with the
estimates coming from the Levi method for this case.

\begin{lem}
Let $n=1$. Suppose that the function $K(t,x;\tau ,\xi )$ is
continuous on
$$
\mathcal P_1=\left\{ (t,x;\tau ,\xi ):\ x,\xi \in \mathbb R^1,x\ne \xi,
0\le \tau <t\le T\right\},
$$
and satisfies the inequality
$$
|K(t,x;\tau ,\xi )|\le A_1(t-\tau )^{\gamma \beta -\beta -1}
\exp \{-\mu \rho (t,x;\tau ,\xi )\}
$$
where $A_1,\mu>0$, $0<\gamma,\beta <1$.Then the series
$\mathcal R(t,x;\tau ,\xi )=\sum\limits_{m=1}^\infty K_m(t,x;\tau ,\xi
)$ converges absolutely and uniformly on $\mathcal P_1\cap
\{(t,x;\tau ,\xi ):\ t-\tau \ge \delta >0\}$, and
\begin{equation}
|\mathcal R(t,x;\tau ,\xi )|\le C(t-\tau )^{\gamma \beta -\beta -1}
\exp \{-\mu^*\rho (t,x;\tau ,\xi )\},
\end{equation}
with some $\mu^*\in (0,\mu )$.
\end{lem}

{\it Proof}. Using Lemma 2, we find for any $\varepsilon >0$ that
\begin{multline*}
|K_2(t,x;\tau ,\xi)|\le
A_1^2\int\limits_\tau^t[(t-\lambda)(\lambda -\tau)]^{\gamma \beta
-\beta-1}\,d\lambda \int\limits_{-\infty}^\infty \exp \{\mu [\rho
(t,x;\lambda ,y)+\rho (\lambda ,y;\tau ,\xi )]\,dy\\
\le A_2(\varepsilon )(t-\tau )^{2\gamma \beta -\beta-1}
\exp \{-\mu (1-\varepsilon )\rho (t,x;\tau ,\xi )\}.
\end{multline*}
Repeating the procedure we obtain the inequality
$$
|K_m(t,x;\tau ,\xi)|\le A_m(\varepsilon )(t-\tau )^{m\gamma \beta -\beta-1}
\exp \{-\mu (1-(m-1)\varepsilon )\rho (t,x;\tau ,\xi )\}.
$$

Let $m^*=\left[ \dfrac{1+\beta }{\beta \gamma }\right] +1$. Then
\begin{equation}
|K_{m^*}(t,x;\tau ,\xi)|\le A_{m^*}(\varepsilon )
\exp \{-\mu (1-m^*\varepsilon )\rho (t,x;\tau ,\xi )\}.
\end{equation}
Setting $\varepsilon =\dfrac{1}{pm^*}$, $p>2$, $\mu^*=\mu
(1-p^{-1})$, $A_0=\max\limits_{1\le m\le m^*}A_m(\varepsilon )$,
we can rewrite (3.32) as
\begin{equation}
|K_{m^*}(t,x;\tau ,\xi)|\le A_0\exp \{-\mu^*\rho (t,x;\tau ,\xi )\}.
\end{equation}

Now, for $m>m^*$, we proceed in a different way, in order to
preserve the exponential factor $\mu^*$ in the estimates of
further iterated kernels. By Lemma 1,
\begin{multline*}
|K_{m^*+1}(t,x;\tau ,\xi)|\le
A_0^2\int\limits_\tau^t(t-\lambda)^{\gamma \beta
-\beta-1}\,d\lambda \int\limits_{-\infty}^\infty \exp \{\mu [\rho
(t,x;\lambda ,y)+(1-p^{-1})\rho (\lambda ,y;\tau ,\xi )]\,dy\\
\le A_0^2F\exp \{-\mu^*\rho (t,x;\tau ,\xi )\}
\int\limits_\tau^t(t-\lambda)^{\gamma \beta -1}\,d\lambda
=A_0^2FB(\gamma \beta ,1)(t-\tau)^{\gamma \beta }\exp
\{-\mu^*\rho (t,x;\tau ,\xi )\},
\end{multline*}
where $F=\int\limits_0^\infty \exp \left\{
-\dfrac{\mu}p|y|^{\frac{1}{1-\beta}}\right\}\,dy$.

Similarly, we prove by induction that
$$
|K_{m^*+j}(t,x;\tau ,\xi)|\le A_0(A_0F)^j(t-\tau )^{\gamma \beta
j}\left\{ \prod\limits_{s=0}^{j-1}B(\gamma \beta ,1+s\gamma \beta
)\right\}\exp \{-\mu^*\rho (t,x;\tau ,\xi )\}.
$$
Using the identity
$$
\prod\limits_{s=0}^{j-1}B(\gamma \beta ,1+s\gamma \beta
)=\frac{[\Gamma (\gamma \beta )]^j}{\Gamma (1+\gamma \beta j)}
$$
we see that the series for $\mathcal R$ is majorized by the convergent
series
$$
\left\{ \sum\limits_{m=1}^{m^*}A_m(t-\tau )^{\gamma m\beta -\beta
-1}+A_0\sum\limits_{j=1}^\infty \frac{[\Gamma (\gamma \beta
)A_0F]^j}{\Gamma (1+\gamma \beta j)}(t-\tau )^{\gamma \beta
j}\right\}\exp \{-\mu^*\rho (t,x;\tau ,\xi )\}.
$$
This implies (3.31). $\quad \Box$

\section{PARAMETRIX}

{\bf 4.1.} {\it The function $Z_0$}. The parametrix kernel
$Z_0(t,x-\xi ,\zeta )$ is defined by the formula (2.3) where
$\left( A^{(ij)}\right)$ is the matrix inverse to the matrix
$(a_{ij}(\zeta ))$ of the leading coefficients ``frozen'' at the
parametric point $\zeta \in \mathbb R^n$. The estimates of $Z_0$
proved in \cite{K2} for the case of the constant coefficients
$a_{ij}$ and collected in Proposition 1, remain valid for the
parametrix, with the constants independent on $\zeta$.

Since properties of $Z_0$ are different for $n\ge 2$ and $n=1$, it
is convenient to treat these cases separately. Therefore we assume
in this and the next sections that $n\ge 2$. The case $n=1$ will
be considered in Sect. 6.

The behavior of $Z_0$ for $R\ge 1$ and $R<1$ is described in
Proposition 1 separately. However it is possible to write
equivalent unified estimates:
\begin{equation}
\left| D_x^mZ_0(t,x-\xi ;\zeta )\right| \le Ct^{-\alpha}|x-\xi
|^{-n+2-|m|}\exp \{-\sigma \rho (t,x;0,\xi )\},\quad |m|\le 3,
\end{equation}
\begin{equation}
\left| \D Z_0(t,x-\xi ;\zeta )\right| \le Ct^{-\alpha}|x-\xi
|^{-n}\exp \{-\sigma \rho (t,x;0,\xi )\},\quad |m|\le 3,
\end{equation}
where
$$
\rho (t,x;\tau ,\xi )=\left( \frac{|x-\xi |}{(t-\tau )^{\alpha
/2}}\right)^{\frac{2}{2-\alpha }}
$$
(this corresponds to the notation of Sect. 3.2 with $\beta =\alpha
/2$). If $n=2$, then (4.1) is valid only for $|m|\ne 0$. For the
opposite case
\begin{equation}
\left| Z_0(t,x-\xi ;\zeta )\right| \le Ct^{-\alpha}\left[ \left| \log
\left( t^{-\alpha }|x-\xi |^2\right) \right|+1\right]
\exp \{-\sigma \rho (t,x;0,\xi )\}.
\end{equation}
The constants $C,\sigma >0$ in (4.1)-(4.3) do not depend on
$\zeta$. Of course, the constants may be different from the ones
in (2.5)-(2.10).

The next proposition gives similar estimates for differences of
values of $Z_0$ and its derivatives corresponding to different
values of the parameter $\zeta$.

\begin{prop}
For any $y,\zeta',\zeta''\in \mathbb R^n$, $0<t\le T$,
\begin{multline}
\left| D_x^mZ_0(t,y;\zeta')-D_x^mZ_0(t,y;\zeta'')\right| \\
\le Ct^{-\alpha}|\zeta'-\zeta''|^\gamma |y|^{-n+2-|m|}\exp
\{-\sigma \left( t^{-\alpha
/2}|y|\right)^{\frac{2}{2-\alpha }}\},\quad |m|\le 2.
\end{multline}
\end{prop}

{\it Proof}. Let $m=0$. Denote
$$
\mathfrak A(y,\zeta )=\sum\limits_{i,j=1}^nA^{(ij)}(\zeta )y_iy_j.
$$
By our assumptions
\begin{gather*}
C_1|y|^2\le \mathfrak A(y,\zeta )\le C_2|y|^2,\\
\left| \mathfrak A(y,\zeta')-\mathfrak A(y,\zeta'')\right| \le
C|\zeta'-\zeta''|^\gamma |y|^2,\\
\left| \left[ \det A(\zeta')\right]^{1/2}-\left[ \det
A(\zeta'')\right]^{1/2}\right| \le C|\zeta'-\zeta''|^\gamma.
\end{gather*}

According to (2.3), we have to use the estimate (3.8) for the
function $H_{12}^{20}$, and also to find an estimate for a
difference of such functions. Let
$$
\varphi (s)=s^{-n/2}H_{12}^{20}\left[ \omega s\left|
\begin{matrix} (1,\alpha ) & \\
(\frac{n}2,1), & (1,1)\end{matrix}\right. \right],\quad \omega
=\frac{1}4t^{-\alpha }.
$$
By (3.1),
$$
\varphi'(s)=-\frac{n}2s^{-\frac{n}2-1}H_{12}^{20}\left[ \omega s\left|
\begin{matrix} (1,\alpha ) & \\
(\frac{n}2,1), & (1,1)\end{matrix}\right. \right] -
s^{-\frac{n}2-1}H_{23}^{30}\left[ \omega s\left|
\begin{matrix} (1,\alpha ), & (0,1) & \\
(\frac{n}2,1), & (1,1), & (1,1)\end{matrix}\right. \right].
$$

For large $s$, (3.7) and (3.8) give the estimate
$$
|\varphi'(s)|\le Cs^{\frac{n+2}{2(2-\alpha )}-\frac{n}2-1}
\omega^{\frac{n+2}{2(2-\alpha )}}\exp \left( -\sigma |\omega
s|^{\frac{1}{2-\alpha }}\right),
$$
which implies the bound
$$
\left| \varphi \left( \mathfrak A(y,\zeta')\right)-\varphi \left(
\mathfrak A(y,\zeta'')\right) \right| \le Ct^{-\frac{n\alpha }2}
|\zeta'-\zeta''|^\gamma \exp \left( -\sigma t^{-\frac{\alpha
}{2-\alpha }}|y|^{\frac{2}{2-\alpha }}\right).
$$

For small $s$, we have
$$
\left| \varphi'(s)\right| \le C\omega s^{-\frac{n}2}
$$
by (3.10) and (3.11); note that the logarithmic terms emerging for
$n=2$ are cancelled.

As a result, we find that
$$
\left| Z_0(t,x-\xi ;\zeta')-Z_0(t,x-\xi ,\zeta'')\right| \le
\begin{cases}
C|\zeta'-\zeta''|^\gamma t^{-\frac{n\alpha}2}\exp \left(
-\sigma t^{-\frac{\alpha}{2-\alpha}}|y|^{\frac{2}{2-\alpha}}\right) &
\text{ if $t^{-\alpha}|y|\ge 1$},\\
C|\zeta'-\zeta''|^\gamma t^{-\alpha}|y|^{-n+2}, &
\text{ if $t^{-\alpha}|y|\le 1$},\end{cases}
$$
which is equivalent to the estimate (4.4) with $m=0$.

For the first and second derivatives the proof is similar, though
somewhat cumbersome -- one has to use the asymptotics (3.7)-(3.12)
for the H-functions; again for $n=2$ the logarithmic terms are
cancelled. $\quad \Box$

\medskip
{\bf 4.2.} {\it The function $Y_0$}. The estimates for the function
$Y_0$ are given in Proposition 2 for the case of constant
coefficients. They carry over to the kernel $Y_0(t-\lambda
,x-\xi;\zeta )$ defined by (2.4) with the coefficients $A^{(ij)}$
depending on $\zeta$ as above. The proof is a direct consequence
of the formula (2.4) and properties of H-functions, in particular
the differentiation formulas (3.1), (3.2) etc, the asymptotic
relations (3.7) and (3.14)-(3.17). The calculations are simple but
tedious, especially for $n\le 4$, since we have to take into
account several terms of the asymptotics of H-functions near the
origin.

Just as for the function $Z_0$, it is desirable to obtain unified
estimates of $Y_0$ and its derivatives which are valid for all
values of the independent variables. Here the problem is a little
more complicated because the behavior of $Y_0$ is different for
different values of the dimension $n$.

In particular, for $n=2$ we have
\begin{equation}
\left| D_x^mY_0(t-\tau ,x-\xi ,\zeta )\right|
\le C(t-\tau )^{-\frac{\alpha}2(2+|m|)+\alpha -1}
\exp \{-\sigma \rho (t,x;\tau ,\xi )\},\quad |m|=0,1;
\end{equation}
\begin{multline}
\left| D_x^mY_0(t-\tau ,x-\xi ;\zeta )\right| \\
\le C(t-\tau )^{-\alpha -1}|x-\xi |^{2-|m|}\left[ \left| \log
\frac{|x-\xi |^2}{(t-\tau )^\alpha}\right| +1\right]
\exp \{-\sigma \rho (t,x;\tau ,\xi )\},\quad |m|=2,3.
\end{multline}
These estimates follow immediately from (2.11), (2.21)-(2.24).

For $n=3$,
\begin{equation}
\left| D_x^mY_0(t-\tau ,x-\xi ;\zeta )\right|
\le C(t-\tau )^{-\frac{\alpha}2(3+|m|)+\alpha -1}
\exp \{-\sigma \rho (t,x;\tau ,\xi )\},\quad |m|=0,1,
\end{equation}
and
\begin{equation}
\left| D_x^mY_0(t-\tau ,x-\xi ;\zeta )\right|
\le C(t-\tau )^{-\alpha -1}|x-\xi |^{1-|m|}
\exp \{-\sigma \rho (t,x;\tau ,\xi )\},\quad |m|=2,3.
\end{equation}

Similarly, for $n=4$,
\begin{multline}
\left| D_x^mY_0(t-\tau ,x-\xi ;\zeta )\right| \\
\le C(t-\tau )^{-\alpha -\frac{|m|\alpha }2-1}\left[ \left| \log
\frac{|x-\xi |^2}{(t-\tau )^\alpha}\right| +1\right]
\exp \{-\sigma \rho (t,x;\tau ,\xi )\},\quad |m|=0,1,2,
\end{multline}
\begin{multline}
\left| D_x^mY_0(t-\tau ,x-\xi ;\zeta )\right| \\
\le C(t-\tau )^{-2\alpha -1}|x-\xi |^{-1}\left[ \left| \log
\frac{|x-\xi |^2}{(t-\tau )^\alpha}\right| +1\right]
\exp \{-\sigma \rho (t,x;\tau ,\xi )\},\quad |m|=3.
\end{multline}

Finally, if $n>4$, then
\begin{equation}
\left| D_x^mY_0(t-\tau ,x-\xi ;\zeta )\right|
\le C(t-\tau )^{-\alpha -1}|x-\xi |^{-n-|m|+4}
\exp \{-\sigma \rho (t,x;\tau ,\xi )\},\quad |m|\le 3.
\end{equation}

Just as in Proposition 3, the estimates for the differences
$D_x^mY_0(t-\tau ,x-\xi ;\zeta')-D_x^mY_0(t-\tau ,x-\xi ;\zeta'')$ have
the following form: the right-hand sides of (4.5)-(4.11) are multiplied
by $|\zeta'-\zeta''|^\gamma$. For example, if $n>4$, then
\begin{multline}
\left| D_x^mY_0(t-\tau ,x-\xi ;\zeta')-D_x^mY_0(t-\tau ,x-\xi
;\zeta'')\right| \\
\le C(t-\tau )^{-\alpha -1}
|\zeta'-\zeta''|^\gamma|x-\xi |^{-n-|m|+4}
\exp \{-\sigma \rho (t,x;\tau ,\xi )\},\quad |m|\le 3.
\end{multline}

{\bf 4.3.} {\it Integral identities}. It follows from the
construction of the function $Z_0$ that
\begin{equation}
\int\limits_{\mathbb R^n}Z_0(t,x-\xi ;\zeta )\,d\xi =1
\end{equation}
(see \cite{K2}). Next, from (2.4) we see that
\begin{multline*}
\int\limits_{\mathbb R^n}Y_0(t-\tau ,x-\xi ;\zeta )\,d\xi =(t-\tau
)^{\alpha -1}\pi^{-n/2}\int\limits_{\mathbb R^n}|y|^{-n}
H_{12}^{20}\left[ \frac{1}{4}(t-\tau )^{-\alpha}|y|^2
\left| \begin{matrix} (\alpha ,\alpha ) & \\
(\frac{n}2,1), & (1,1)\end{matrix}\right. \right]\,dy\\
=C(t-\tau )^{\alpha -1}\int\limits_0^\infty r^{-1}
H_{12}^{20}\left[ \frac{1}{4}(t-\tau )^{-\alpha}r^2
\left| \begin{matrix} (\alpha ,\alpha ) & \\
(\frac{n}2,1), & (1,1)\end{matrix}\right. \right]\,dr\\ =C_1
(t-\tau )^{\alpha -1}\int\limits_0^\infty s^{-1}
H_{12}^{20}\left[ \frac{1}{4}(t-\tau )^{-\alpha}s
\left| \begin{matrix} (\alpha ,\alpha ) & \\
(\frac{n}2,1), & (1,1)\end{matrix}\right. \right]\,ds,
\end{multline*}
so that
\begin{equation}
\int\limits_{\mathbb R^n}Y_0(t-\tau ,x-\xi ,\zeta )\,d\xi =C_2(t-\tau)^{\alpha -1}.
\end{equation}

\medskip
{\bf 4.4.} {\it Further estimates}. Below we will need also the
first time derivative of the function $Z_0(t,x;\zeta )$. We find
from (2.3) and (3.1) that
\begin{multline*}
\frac{\partial Z_0(t,x;\zeta )}{\partial t}=\frac{\alpha
\pi^{-n/2}t^{-1}}{(\det A)^{1/2}}\left[
\sum\limits_{i,j=1}^nA^{(ij)}x_ix_j\right]^{-n/2}\\
\times H_{23}^{30}\left[ \frac{1}{4}t^{-\alpha}
\sum\limits_{i,j=1}^nA^{(ij)}x_ix_j\left|
\begin{matrix} (1,\alpha ), & (0,1) & \\
(\frac{n}2,1), & (1,1), & (1,1)\end{matrix}\right. \right].
\end{multline*}
Using the asymptotics (3.7) and (3.11) of the function
$H_{23}^{30}$ we can obtain an estimate of $\dfrac{\partial
Z_0}{\partial t}$. For example, if $n\ge 3$, then
\begin{equation}
\left| \frac{\partial Z_0(t,x;\zeta )}{\partial t}\right| \le
\begin{cases}
Ct^{-1-\frac{\alpha n}2}\exp \left\{-\sigma (t^{-\alpha
/2}|x|)^{\frac{2}{2-\alpha }}\right\}, &
\text{ if $t^{-\alpha /2}|x|\ge 1$},\\
Ct^{-\alpha -1}|x|^{-n+2}, & \text{ if $t^{-\alpha /2}|x|<1$}.
\end{cases}
\end{equation}

We can also use a roughened unified estimate
$$
\left| \frac{\partial Z_0(t,x;\zeta )}{\partial t}\right| \le
Ct^{-1}|x|^{-n}\exp \left\{-\sigma (t^{-\alpha
/2}|x|)^{\frac{2}{2-\alpha }}\right\}.\eqno (4.15')
$$

As before, we get also an estimate for the difference
$\left|\dfrac{\partial Z_0(t,x;\zeta')}{\partial t}-\dfrac{\partial
Z_0(t,x;\zeta'')}{\partial t}\right|$ whose upper bound is the
expression in the right-hand side of (4.15) or $(4.15')$
multiplied by $|\zeta'-\zeta''|^\gamma$. Since
$$
\int\limits_{\mathbb R^n}\frac{\partial Z_0(t,x-\xi ;\zeta')}{\partial
t}\,d\xi =0
$$
by virtue of (4.13), the above estimates imply also the estimate
\begin{equation}
\left| \int\limits_{\mathbb R^n}\frac{\partial Z_0(t,x-\xi ;\xi )}{\partial
t}\,d\xi \right| \le Ct^{-1+\frac{\alpha \gamma }2}
\end{equation}
obtained by subtracting the (zero) integral of
$\dfrac{\partial Z_0(t,x;x)}{\partial t}$ and using the estimate
for the difference of the derivatives.

\medskip
\section{THE LEVI METHOD ($n\ge 2$)}

{\bf 5.1.} {\it The scheme}. We look for the functions $Y,Z$ appearing in
(1.4) assuming the following integral representations:
\begin{equation}
Z(t,x;\xi )=Z_0(t,x-\xi ;\xi )+\int\limits_0^td\lambda
\int\limits_{\mathbb R^n}Y_0(t-\lambda ,x-y;y)Q(\lambda ,y;\xi
)\,dy;
\end{equation}
\begin{equation}
Y(t,x;\xi )=Y_0(t,x-\xi ;\xi )+\int\limits_0^td\lambda
\int\limits_{\mathbb R^n}Y_0(t-\lambda ,x-y;y)\Psi (\lambda ,y;\xi
)\,dy;
\end{equation}
the functions $Z_0,Y_0$ were examined in detail in Sect. 4.

For the functions $Q$ and $\Psi$ we assume the integral equations
\begin{equation}
Q(t,x;\xi )=M(t,x;\xi )+\int\limits_0^td\lambda
\int\limits_{\mathbb R^n}K(t-\lambda,x;y)Q(\lambda ,y;\xi
)\,dy,
\end{equation}
\begin{equation}
\Psi (t,x;\xi )=K(t,x;\xi )+\int\limits_0^td\lambda
\int\limits_{\mathbb R^n}K(t-\lambda,x;y)\Psi (\lambda ,y;\xi )\,dy,
\end{equation}
where
\begin{multline*}
M(t,x;\xi )=\sum\limits_{i,j=1}^n\left\{ [a_{ij}(x)-a_{ij}(\xi
)]\frac{\partial^2}{\partial x_i\partial x_j}Z_0(t,x-\xi ;\xi
)\right\}\\
+\sum\limits_{i=1}^nb_i(x)\frac{\partial Z_0(t,x-\xi ;\xi
)}{\partial x_i}+c(x)Z_0(t,x-\xi ;\xi ),
\end{multline*}
\begin{multline*}
K(t,x;\xi )=\sum\limits_{i,j=1}^n\left\{ [a_{ij}(x)-a_{ij}(\xi
)]\frac{\partial^2}{\partial x_i\partial x_j}Y_0(t,x-\xi ;\xi
)\right\}\\
+\sum\limits_{i=1}^nb_i(x)\frac{\partial Y_0(t,x-\xi ;\xi
)}{\partial x_i}+c(x)Y_0(t,x-\xi ;\xi).
\end{multline*}

Using the estimates (4.1)-(4.3) we find that
\begin{equation}
|M(t,x;\xi )|\le Ct^{-\alpha }|x-\xi |^{-n+\gamma }\exp \{-\sigma
\rho (t,x;0,\xi )\}.
\end{equation}

In order to obtain estimates for $K$, we have to use estimates for
$Y_0$, different for different values of $n$, and for different
domains, and then to roughen the resulting estimates into less
exact but unified bounds convenient for the Levi method.

\begin{prop}
For any $n\ge 2$
\begin{equation}
|K(t,x;\xi )|\le Ct^{(\gamma -\eta )\frac{\alpha}2-1}
|x-\xi |^{-n+\eta }\exp \{-\sigma \rho (t,x;0,\xi )\},
\end{equation}
with $0<\eta <\gamma$.
\end{prop}

{\it Proof}. Using Proposition 2 we write the estimate for $|x-\xi
|\ge t^{\alpha /2}$:
\begin{equation}
|K(t,x;\xi )|\le Ct^{-\frac{\alpha n}2+\frac{\alpha
\gamma }2-1}\exp \{-\sigma \rho (t,x;0,\xi )\}.
\end{equation}

For $|x-\xi |\le t^{\alpha /2}$ we consider various cases
separately. If $n=3$, or $n>4$, then
\begin{equation}
|K(t,x;\xi )|\le Ct^{-\alpha +\frac{\alpha \gamma}2-1}
|x-\xi |^{-n+2}\exp \{-\sigma \rho (t,x;0,\xi )\}.
\end{equation}
Since
$$
|x-\xi |^2=|x-\xi|^\eta |x-\xi|^{2-\eta}\le |x-\xi|^\eta
t^{\frac{\alpha}2(2-\eta)},
$$
(5.7) and (5.8) imply (5.6). Note that transforming (5.7) we
change the constant $\sigma$.

If $n=2$ or $n=4$, the initial bound for $K$ is
\begin{equation}
|K(t,x;\xi )|\le Ct^{-\frac{\alpha n}2+\frac{\alpha
\gamma}2-1}\left( \left| \log \frac{|x-\xi |^2}{t^\alpha}
\right|+1\right) \exp \{-\sigma \rho (t,x;0,\xi )\}.
\end{equation}
We roughen the estimate (5.9) replacing the factor
$\left| \log \dfrac{|x-\xi |^2}{t^\alpha}\right|+1$
with the factor \linebreak
$\left( \dfrac{|x-\xi |}{t^{\alpha/2}}\right)^{-n+\eta}$,
which results in (5.6).$\quad \Box$

\medskip
Next we have to study the increments
\begin{gather*}
\Delta_xM(t,x;\xi )=M(t,x;\xi )-M(t,x';\xi ),\\
\Delta_xK(t,x;\xi )=K(t,x;\xi )-K(t,x';\xi ),
\end{gather*}
$x,x'\in \mathbb R^n$.

We can write $\Delta_xM(t,x;\xi )=M_1+M_2+M_3$ where
$$
M_1=\sum\limits_{i,j=1}^n\Delta_xa_{ij}(x)
\frac{\partial^2Z_0(t,x-\xi ;\xi)}{\partial x_i\partial x_j},
$$
$$
M_2=\sum\limits_{i,j=1}^n[a_{ij}(x')-a_{ij}(\xi)]\Delta_x
\frac{\partial^2Z_0(t,x'-\xi ;\xi )}{\partial x_i'\partial x_j'},
$$
and the term $M_3$ contains lower order derivatives. Since $M_3$
has a weaker singularity and does not influence the estimates, we
omit its detailed description and consider only $M_1$ and $M_2$.

An estimate for $M_1$ is given directly:
\begin{equation}
|M_1|\le C_1|x-x'|^{\gamma}t^{-\alpha }|x-\xi |^{-n}\exp \{-\sigma
\rho (t,x;0,\xi )\}.
\end{equation}

Using (4.1) we find that
\begin{multline}
\left| \Delta_x
\frac{\partial^2Z_0(t,x'-\xi ;\xi )}{\partial x_i'\partial
x_j'}\right|\le \left| \Delta_x\frac{\partial^2Z_0}{\partial x_i'\partial
x_j'}\right|^\nu \left( \left| \frac{\partial^2Z_0}{\partial x_i\partial
x_j}\right| +\left| \frac{\partial^2Z_0}{\partial x_i'\partial
x_j'}\right|\right)^{1-\nu}\\
\le C_3^\nu t^{-\alpha \nu}|\widetilde{x}-\xi |^{-n\nu -\nu}
\exp \{-\nu \sigma \rho (t,\widetilde{x};0,\xi )\}|x-x'|^\nu
C_2^{1-\nu}t^{-\alpha (1-\nu )}\left[ |x-\xi |^{-n}\exp \{-\sigma \rho
(t,x;0,\xi )\}\right. \\
\left. +|x'-\xi |^{-n}\exp \{-\sigma \rho (t,x';0,\xi
)\}\right]^{1-\nu}
\end{multline}
where $\widetilde{x}=x+\theta (x'-x)$, $0<\theta <1$, $0<\nu
<\gamma$, and we assumed that $\min \left\{ |x-\xi |,|x'-\xi |\right\}
\ne 0$. From (5.10) and (5.11) we find that
\begin{equation}
|\Delta_xM(t,x;\xi )|\le C_4t^{-\alpha }\left[ |x-x'|^\gamma |x-\xi |^{-n}
\exp \{-\sigma \rho (t,x;0,\xi )\} +|x-x'|^\nu |x'-\xi
|^\gamma \widetilde{M}(t,x;\xi )\right]
\end{equation}
where
\begin{multline*}
\widetilde{M}(t,x;\xi )=|\widetilde{x}-\xi |^{-n\nu -\nu}
\exp \{-\nu \sigma \rho (t,\widetilde{x};0,\xi )\}
 \left[ |x-\xi |^{-n}\exp \{-\sigma \rho
(t,x;0,\xi )\}\right. \\
\left. +|x'-\xi |^{-n}\exp \{-\sigma \rho (t,x';0,\xi
)\}\right]^{1-\nu}.
\end{multline*}
We will use also the following estimate which is a direct
consequence of the definition of $\Delta_xM$:
\begin{multline}
|\Delta_xM(t,x;\xi )|\le Ct^{-\alpha }\left( |x-\xi |^{-n+\gamma}
\exp \{-\sigma \rho (t,x;0,\xi )\}\right. \\
\left. +|x'-\xi |^{-n+\gamma}
\exp \{-\sigma \rho (t,x';0,\xi )\}\right) .
\end{multline}

Let us consider two possible cases.

a) Suppose that $|x-x'|>At^{\alpha /2}$. It follows from (5.13)
that
\begin{multline*}
|\Delta_xM(t,x;\xi )|\le C_5(A)t^{(\gamma-\varepsilon)\frac{\alpha}2-\alpha
}\left[ |x-\xi |^{-n+\varepsilon }\exp \{-\sigma_1\rho (t,x;0,\xi
)\}\right. \\
\left. +|x'-\xi |^{-n+\varepsilon }
\exp \{-\sigma_1\rho (t,x';0,\xi )\}\right],\quad \varepsilon >0,
\end{multline*}
whence
\begin{multline}
|\Delta_xM(t,x;\xi )|\le C_6t^{-\alpha }|x-x'|^{\gamma -\varepsilon }
\left[ |x-\xi |^{-n+\varepsilon }\exp \{-\sigma_1\rho (t,x;0,\xi
)\}\right. \\
\left. +|x'-\xi |^{-n+\varepsilon }
\exp \{-\sigma_1\rho (t,x';0,\xi )\}\right],\quad \varepsilon >0,
\end{multline}

b) Let $|x-x'|\le At^{\alpha /2}$ with some $A>0$. If $|x-x'|>\eta
\max (|\xi -x|,|\xi -x'|)$ with some $\eta >0$, then we obtain
again the estimate (5.14). Thus, we have now to consider the most
complicated subcase, for which
$$
|x-x'|\le \eta \max (|\xi -x|,|\xi -x'|),\quad \eta >0.
$$
Here we use the estimate (5.11).

Without restricting generality, we may assume that $|\xi -x'|\le
|\xi -x|$. Then $|x-x'|\le \eta |\xi -x|$,
$$
\left| \widetilde{x}-\xi \right|=|x-\xi +\theta (x'-x)|\ge |x-\xi
|-\theta \eta |x-\xi |\ge (1-\eta )|x-\xi |.
$$
Note that the function
$$
r\mapsto r^{-n}\exp \left\{ -\sigma \left( \frac{r}{t^{\alpha
/2}}\right)^{\frac{2}{2-\alpha }}\right\} ,\quad 0<\alpha <1,
$$
is monotone decreasing. This implies the estimate
$$
\left| \widetilde{M}(t,x;\xi )\right| \le C_7\left[
|x-\xi |^{-n-\nu }\exp \{-\sigma_1 \rho (t,x;0,\xi
)\} +|x'-\xi |^{-n-\nu }\exp \{-\sigma \rho (t,x';0,\xi )\}\right] .
$$
Setting $\varepsilon =\gamma -\nu$ and using (5.13), (5.14), we
come to the estimate contained in the following proposition.

\begin{prop}
If $x''$ is one of the points $x,x'$ for which $|x''-\xi |=\min \{
|x-\xi |,|x'-\xi |\}$, then
\begin{equation}
|\Delta_xM(t,x;\xi )|\le Ct^{-\alpha }|x-x'|^{\gamma -\varepsilon
}|x''-\xi |^{-n+\varepsilon }\exp \{-\sigma \rho (t,x'';0,\xi )\},
\end{equation}
\begin{equation}
|\Delta_xK(t,x;\xi )|\le Ct^{-1}|x-x'|^{\gamma -\varepsilon
}|x''-\xi |^{-n+\varepsilon }\exp \{-\sigma \rho (t,x'';0,\xi
)\}.
\end{equation}
\end{prop}

{\it The proof} of (5.16) is similar to that of (5.15), since the estimates
for $Y_0$ have a structure similar to those for $Z_0$. $\quad \Box$

\medskip
Now we are ready to consider the integral equations (5.3) and (5.4).

\begin{prop}
The integral equations (5.3) and (5.4) have the solutions $Q(t,x;\xi
)$, $\Psi (t,x;\xi )$, which are continuous for
$x\ne \xi$ and satisfy the estimates
\begin{equation}
|Q(t,x;\xi )|\le Ct^{-\alpha }|x-\xi |^{-n+\gamma }\exp
\{-\sigma \rho (t,x;0,\xi )\},
\end{equation}
\begin{equation}
|\Delta_xQ(t,x;\xi )|\le Ct^{-\alpha }|x-x'|^{\gamma -\varepsilon
}|x''-\xi |^{-n+\varepsilon }\exp \{-\sigma \rho (t,x'';0,\xi )\},
\end{equation}
\begin{equation}
|\Psi (t,x;\xi )|\le Ct^{-1}|x-\xi |^{-n+\gamma }\exp
\{-\sigma \rho (t,x;0,\xi )\},
\end{equation}
\begin{equation}
|\Delta_x\Psi (t,x;\xi )|\le Ct^{-1}|x-x'|^{\gamma -\varepsilon
}|x''-\xi |^{-n+\varepsilon }\exp \{-\sigma \rho (t,x'';0,\xi )\},
\end{equation}
with an arbitrary $\varepsilon \in (0,\gamma )$ and the constants
depending on $\varepsilon$. If $n=2$ or $n=4$, then for any $\mu >0$
\begin{equation}
|Q(t,x;\xi )|\le Ct^{-\alpha -\frac{\mu \alpha}2}
|x-\xi |^{-n+\varepsilon -\mu }\exp \{-\sigma \rho (t,x;0,\xi )\},
\tag{$5.17'$}
\end{equation}
\begin{equation}
|\Delta_xQ(t,x;\xi )|\le Ct^{-\alpha -\frac{\mu \alpha}2
}|x-x'|^{\gamma -\varepsilon
}|x''-\xi |^{-n+\varepsilon -\mu }\exp \{-\sigma \rho (t,x'';0,\xi )\},
\tag{$5.18'$}
\end{equation}
\begin{equation}
|\Psi (t,x;\xi )|\le Ct^{-1-\frac{\mu \alpha}2}
|x-\xi |^{-n+\varepsilon -\mu }\exp \{-\sigma \rho (t,x;0,\xi )\},
\tag{$5.19'$}
\end{equation}
\begin{equation}
|\Delta_x\Psi (t,x;\xi )|\le Ct^{-1-\frac{\mu \alpha}2}
|x-x'|^{\gamma -\varepsilon
}|x''-\xi |^{-n+\varepsilon -\mu }\exp \{-\sigma \rho (t,x'';0,\xi
)\}.
\tag{$5.20'$}
\end{equation}
Here the constants depend on $\varepsilon$ and $\mu$.
\end{prop}

{\it Proof}. The equations (5.3) and (5.4) are solved by
iterations. The estimates $(5.17)$, $(5.17')$, $(5.19)$, $(5.19')$
follow from Lemma 5.

Let us prove (5.18). We have $\Delta_xQ=\Delta_xM+J$ where
\begin{equation}
J=\int\limits_0^td\lambda \int\limits_{\mathbb
R^n}\Delta_xK(t-\lambda,x;y)Q(\lambda ,y;\xi )\,dy,
\end{equation}
\begin{multline*}
|J|\le C|x-x'|^{\gamma -\varepsilon }\int\limits_0^t(t-\lambda )^{-1}
\lambda^{-\alpha }d\lambda \left[ \int\limits_{\mathbb
R^n}|x-y|^{-n+\varepsilon }|y-\xi |^{-n+\gamma }\exp \{ -\sigma
\rho (t,x;\lambda ,y)\right. \\
-\sigma \rho (\lambda ,y;0,\xi )\}\,dy
\left. +\int\limits_{\mathbb
R^n}|x'-y|^{-n+\varepsilon }|y-\xi |^{-n+\gamma }\exp \{ -\sigma
\rho (t,x';\lambda ,y)-\sigma \rho (\lambda ,y;0,\xi
)\}\,dy\right] \\
\le C|x-x'|^{\gamma -\varepsilon }\int\limits_0^t(t-\lambda )^{-1+
\frac{\alpha \varepsilon}4}\lambda^{(1-\alpha)-1}\,d\lambda
\left[ \int\limits_{\mathbb
R^n}|x-y|^{-n+\frac{\varepsilon}2}|y-\xi |^{-n+\gamma }\exp \{ -\sigma_1
\rho (t,x;\lambda ,y)\right.
\\  -\sigma_1 \rho (\lambda ,y;0,\xi )\}\,dy
\left. +\int\limits_{\mathbb
R^n}|x'-y|^{-n+\frac{\varepsilon}2}|y-\xi |^{-n+\gamma }\exp \{ -\sigma_1
\rho (t,x';\lambda ,y)-\sigma_1\rho (\lambda ,y;0,\xi
)\}\,dy\right].
\end{multline*}

By Lemma 3, the integral (5.21) is estimated as follows:
\begin{multline*}
|J|\le C|x-x'|^{\gamma -\varepsilon }t^{-\alpha +\frac{\alpha
\varepsilon}4}\left[ |x-\xi |^{-n+\gamma +\frac{\varepsilon}2}\exp
\{ -\sigma_2\rho (t,x;0,\xi )\}\right. \\
\left. +|x'-\xi |^{-n+\gamma +\frac{\varepsilon}2}\exp
\{ -\sigma_2\rho (t,x;0,\xi )\}\right].
\end{multline*}
Together with (5.15), this implies (5.18) with an appropriate
$\sigma >0$.

The proof of the inequalities $(5.20)$, $(5.20')$ is similar. $\quad \Box$

{\bf 5.2.} {\it Heat potentials}. In order to verify that the
functions (5.1) and (5.2) indeed solve the Cauchy problem via the
representation (1.4), we have to study the heat potential
\begin{equation}
W(t,x)=\int\limits_0^td\lambda \int\limits_{\mathbb
R^n}Y_0(t-\lambda ,x-y;y)f(\lambda ,y)\,dy
\end{equation}
for two situations. Firstly, we must consider the case of a
bounded, locally H\"older continuous (in $x$) function $f$, jointly
continuous in $(t,x)\in [0,T]\times \mathbb R^n$. This is, of course,
dictated by the formula (1.4) containing a term which solves the
inhomogeneous equation with the zero initial condition. Note that
$Y_0$ is the main singular part of the function $Y$ appearing in (1.4).
Secondly, heat potentials appear also in (5.1) and (5.2), where
the role of $f$ is played by the functions $Q$ and $\Psi$
which are much more singular (see Proposition 6).

Let us consider the potential (5.22) with a bounded, locally H\"older
continuous $f$. The existence of the integral (5.22) and the possibility
to find its first order derivatives in $x$ by differentiating under
the integral's symbol follow directly from the estimates (4.5), (4.7),
(4.9), and (4.11). For studying the second order derivatives we
consider the function
$$
W(t,x)=\int\limits_0^{t-h}d\lambda \int\limits_{\mathbb
R^n}Y_0(t-\lambda ,x-y;y)f(\lambda ,y)\,dy, \quad h>0.
$$
As $t-\lambda \ge h$, the function $Y_0$ and its first and second
derivatives have integrable singularities at $x=y$. Therefore we
may differentiate under the integral's symbol, so that
$$
\frac{\partial^2W_h(t,x)}{\partial x_i\partial x_j}=J_1+J_2
$$
where
\begin{gather}
J_1=\int\limits_0^{t-h}d\lambda \int\limits_{\mathbb
R^n}\frac{\partial^2Y_0(t-\lambda ,x-y;y)}{\partial x_i\partial x_j}
[f(\lambda ,y)-f(\lambda ,x)]\,dy,\\
J_2=\int\limits_0^{t-h}f(\lambda ,x)\,d\lambda \int\limits_{\mathbb
R^n}\frac{\partial^2Y_0(t-\lambda ,x-y;y)}{\partial x_i\partial
x_j}\,dy.
\end{gather}

We decompose further, $J_1=J_1^{(1)}+J_1^{(2)}$, where $J_1^{(1)}$
corresponds to integration over the set $\Pi^{(1)}=\left\{ y\in
\mathbb R^n:\ |x-y|\le (t-\lambda )^{\alpha /2}\right\}$ while for
$J_1^{(2)}$ the domain of integration is $\Pi^{(2)}=\left\{ y\in
\mathbb R^n:\ |x-y|>(t-\lambda )^{\alpha /2}\right\}$. Suppose that, for
example, $n>4$; other cases are treated similarly. If $y\in \Pi^{(1)}$,
then
$$
\left| \frac{\partial^2Y_0(t-\lambda ,x-y;y)}{\partial x_i\partial
x_j}\right| \le C(t-\lambda )^{-\alpha -1}|x-y|^{-n+2},
$$
so that the $\Pi^{(1)}$ part of the integral in (5.23) is majorized by
$$
I_1^{(1)}=C\int\limits_0^{t-h}(t-\lambda )^{-\alpha -1}d\lambda
\int\limits_{\Pi^{(1)}}|x-y|^{-n+2+\gamma }\,dy,
$$
and the change of variables $\widehat{y}=(t-\lambda )^{-\alpha /2}(y-x)$
gives
\begin{equation}
I_1^{(1)}\le C\int\limits_0^{t-h}(t-\lambda )^{-1+\frac{\alpha
\gamma}2}\,d\lambda \int\limits_{|\widehat{y}|\le 1}
|\widehat{y}|^{-n+2+\gamma }\,d\widehat{y}.
\end{equation}
Here $\gamma >0$ is the H\"older exponent of the function $f(\lambda ,x)$
in $x$.

For $y\in \Pi^{(2)}$ we use the inequality
$$
\left| \frac{\partial^2Y_0(t-\lambda ,x-y;y)}{\partial x_i\partial
x_j}\right| \le C(t-\lambda )^{-\frac{\alpha n}2-1}\exp
\left\{ -\sigma \left(
\frac{|x-y|}{(t-\lambda )^{\alpha /2}}\right)^{\frac{2}{2-\alpha
}}\right\},
$$
a consequence of the general estimate (4.11). A similar argument
shows that the $\Pi^{(2)}$ part of the integral in $J_1$ is
majorized by
$$
C\int\limits_0^{t-h}(t-\lambda )^{-1+\frac{\alpha
\gamma}2}\,d\lambda \int\limits_{|\widehat{y}|>1}
|\widehat{y}|^\gamma \exp \left\{ -\sigma |\widehat{y}|^{\frac{2}{2-\alpha
}}\right\} \,d\widehat{y}.
$$
We see the existence of the limit of $J_1$ as $h\to 0$.

Let us consider $J_2$. Using (4.14) we have
\begin{equation}
J_2=\int\limits_0^{t-h}f(\lambda ,x)\,d\lambda \int\limits_{\mathbb
R^n}\left[ \left.
\frac{\partial^2Y_0(t-\lambda ,x-y;y)}{\partial x_i\partial
x_j}-\frac{\partial^2Y_0(t-\lambda ,x-y;\delta )}{\partial x_i\partial
x_j}\right|_{\delta =x}\right] \,dy.
\end{equation}
By (4.12), the expression in brackets is majorized by
$$
C|x-y|^{-n+2+\gamma }(t-\lambda )^{-\alpha -1}\exp \{-\sigma \rho
(t,x;\lambda ,y)\}.
$$
The same change of variables (used above) shows that the integrand
in (5.26) is majorized by
$$
(t-\lambda )^{-1+\frac{\alpha \gamma}2}|\widehat{y}|^{-n+2+\gamma
}\exp \left\{ -\sigma |\widehat{y}|^{\frac{2}{2-\alpha
}}\right\}.
$$
This means the existence of the limit of $J_2$ as $h\to 0$.
Simultaneously we have proved the formula
\begin{multline}
\frac{\partial^2W(t,x)}{\partial x_i\partial x_j}
=\int\limits_0^td\lambda \int\limits_{\mathbb
R^n}\frac{\partial^2Y_0(t-\lambda ,x-y;y)}{\partial x_i\partial x_j}
[f(\lambda ,y)-f(\lambda ,x)]\,dy\\
+\int\limits_0^tf(\lambda ,x)\,d\lambda \int\limits_{\mathbb
R^n}\frac{\partial^2Y_0(t-\lambda ,x-y;y)}{\partial x_i\partial
x_j}\,dy.
\end{multline}

Now we have to carry out a similar procedure for the potential
(5.22) with $f(\lambda ,y)=Q(\lambda ,y;\xi )$ with a fixed $\xi$
(the case of $\Psi (\lambda ,y;\tau ,\xi )$ is quite similar). We
use the estimates (5.17) and (5.18) for the function $Q$ and its
increment, and also the following roughened estimate:
\begin{equation}
\left| \frac{\partial^2Y_0(t-\lambda ,x-y;\xi )}{\partial x_i\partial
x_j}\right| \le C(t-\lambda )^{-1}|x-y|^{-n}\exp \{-\sigma \rho
(t,x;\lambda ,y)\}.
\end{equation}

In fact we have to prove convergence of the integrals
\begin{gather*}
I'_1=\int\limits_0^td\lambda \int\limits_{\mathbb
R^n}\left| \frac{\partial^2Y_0(t-\lambda ,x-y;y)}{\partial x_i\partial x_j}
\right| |Q(\lambda ,x;\xi )-Q(\lambda ,y;\xi )|\,dy,\\
I'_2=\int\limits_0^t|Q(\lambda ,x;\xi )|\,d\lambda \int\limits_{\mathbb
R^n}\left| \frac{\partial^2Y_0(t-\lambda ,x-y;y)}{\partial x_i\partial
x_j}-\left. \frac{\partial^2Y_0(t-\lambda ,x-y;\delta
)}{\partial x_i\partial x_j}\right|_{\delta =x}\right|\,dy.
\end{gather*}

Let $a>0$ be a small positive constant. Changing $\sigma$ we can
rewrite (5.28) as
$$
\left| \frac{\partial^2Y_0(t-\lambda ,x-y;\xi )}{\partial x_i\partial
x_j}\right| \le C(t-\lambda )^{-1+\frac{a\alpha}2}|x-y|^{-n-a}\exp
\{-\sigma \rho (t,x;\lambda ,y)\}.
$$
Now we get by Lemma 3 and Lemma 4 that for $0<\varepsilon <\gamma$
\begin{multline*}
I_1'\le C\int\limits_0^t(t-\lambda )^{-1+\frac{a\alpha}2}
\lambda^{-\alpha }d\lambda \left[ \int\limits_{\mathbb
R^n}|x-y|^{-n+\gamma -\varepsilon -a}|x-\xi |^{-n+\varepsilon
}\exp \{ -\sigma \rho (t,x;\lambda ,y)\right. \\
\left. -\sigma \rho (\lambda ,x;0,\xi )\}\,dy+\int\limits_{\mathbb
R^n}|x-y|^{-n+\gamma -\varepsilon -a}|y-\xi |^{-n+\varepsilon
}\exp \{ -\sigma \rho (t,x;\lambda ,y)-\sigma \rho (\lambda ,y;0,\xi )\}
\,dy\right]\\
=C|x-\xi |^{-n+\varepsilon }\exp \{-\sigma \rho (t,x;0,\xi )\}
\int\limits_0^t(t-\lambda )^{-1+\frac{a\alpha}2}
\lambda^{-\alpha }d\lambda \\
\times \int\limits_{\mathbb
R^n}|x-y|^{-n+\gamma -\varepsilon -a}\exp \{ -\sigma \rho (t,x;\lambda
,y)\}\,dy\\
+C\int\limits_0^t(t-\lambda )^{-1+\frac{a\alpha}2}
\lambda^{-\alpha }d\lambda \int\limits_{\mathbb
R^n}|x-y|^{-n+\gamma -\varepsilon -a}|y-\xi |^{-n+\varepsilon
}\exp \{ -\sigma \rho (t,x;\lambda ,y)-\sigma \rho (\lambda ,y;0,\xi )\}
\,dy\\
\le C|x-\xi |^{-n+\varepsilon }\exp \{-\sigma \rho (\lambda ,x;0,\xi )\}
\int\limits_0^t(t-\lambda )^{\frac{(\gamma -\varepsilon )\alpha}2-1}
\lambda^{-\alpha }d\lambda \\
+Ct^{-\alpha +\frac{a\alpha}2}|x-\xi
|^{-n+\gamma -a}\exp \{-\sigma \rho (t,x;0,\xi )\}.
\end{multline*}
As before, the letters $\sigma ,C$ meant various positive
constants; this abuse of notation will be convenient in the sequel
too.

In order to obtain a bound for $I_2'$, we use the estimate
\begin{multline*}
\left| \frac{\partial^2Y_0(t-\lambda ,x-y;y)}{\partial x_i\partial
x_j}-\left. \frac{\partial^2Y_0(t-\lambda ,x-y;\delta
)}{\partial x_i\partial x_j}\right|_{\delta =x}\right| \\
\le C|x-y|^{\gamma -n}(t-\lambda )^{-1}\exp \{-\sigma \rho
(t,x;\lambda,y)\}
\le C|x-y|^{-n+\frac{\gamma}2}(t-\lambda )^{-1+\frac{\alpha \gamma}4}
\exp \{-\sigma \rho (t,x;\lambda,y)\}
\end{multline*}
(with a different $\sigma$). Then by Lemma 4
\begin{multline*}
I_2'\le C|x-\xi |^{-n+\gamma}\exp \{-\sigma \rho (t,x;0,\xi )\}\\
\times \int\limits_0^t(t-\lambda )^{-1+\frac{\alpha
\gamma}4}\lambda^{-\alpha}\,d\lambda \int\limits_{\mathbb
R^n}|x-y|^{-n+\frac{\gamma}2}\exp \{-\sigma \rho (t,x;\lambda,y)\}
\\
\le C|x-\xi |^{-n+\gamma}\exp \{-\sigma \rho (t,x;0,\xi )\}
\int\limits_0^t(t-\lambda )^{-1+\frac{\alpha
\gamma}2}\lambda^{-\alpha}\,d\lambda \\
\le Ct^{\frac{\alpha \gamma}2-\alpha}|x-\xi |^{-n+\gamma}\exp
\{-\sigma \rho (t,x;0,\xi )\}.
\end{multline*}

Thus, we have proved (5.27) for $f(\lambda ,y)=Q(\lambda ,y;\xi
)$.

\medskip
{\bf 5.3.} {\it The fractional derivative}. Our next task is to
study the fractional derivative of the heat potential (5.22),
again for the above two cases.

Suppose that $f(t,x)$ is bounded, locally H\"older continuous in $x$, and
jointly continuous in $(t,x)$. We have $\D W=\frac{\partial}{\partial t}v$
where $v=I_{0+}^{1-\alpha}W$ is the Riemann-Liouville fractional
integral, that is
\begin{multline*}
v(t,x)=\frac{1}{\Gamma (1-\alpha
)}\int\limits_0^t(t-\theta)^{-\alpha}\,d\theta
\int\limits_0^\theta d\lambda \int\limits_{\mathbb R^n}Y_0(\theta
-\lambda ,x-y;y)f(\lambda ,y)\,dy\\
=\frac{1}{\Gamma (1-\alpha
)}\int\limits_0^td\lambda \int\limits_\lambda^t
(t-\theta)^{-\alpha}\,d\theta \int\limits_{\mathbb R^n}Y_0(\theta
-\lambda ,x-y;y)f(\lambda ,y)\,dy\\
=\int\limits_0^td\lambda \int\limits_{\mathbb R^n}Z_0(t-\lambda
,x-y;y)f(\lambda ,y)\,dy,
\end{multline*}
by the definition of $Y_0$.

Let $h$ be a small positive number,
$$
v_h(t,x)=\int\limits_0^{t-h}d\lambda \int\limits_{\mathbb R^n}Z_0(t-\lambda
,x-y;y)f(\lambda ,y)\,dy.
$$
Then $v_h\to v$ pointwise as $h\to 0$,
$$
\frac{\partial v_h}{\partial t}=v_h^{(1)}(t,x)+v_h^{(2)}(t,x)
$$
where
\begin{gather*}
v_h^{(1)}(t,x)=\int\limits_{\mathbb
R^n}Z_0(h,x-y;y)f(t-h,y)\,dy,\\
v_h^{(2)}(t,x)=\int\limits_0^{t-h}d\lambda \int\limits_{\mathbb R^n}
\frac{\partial Z_0(t-\lambda ,x-y;y)}{\partial t}f(\lambda ,y)\,dy.
\end{gather*}

We have
\begin{multline}
v_h^{(1)}(t,x)=\int\limits_{\mathbb
R^n}[Z_0(h,x-y;y)-Z_0(h,x-y;x)]f(t-h,y)\,dy\\
+\int\limits_{\mathbb R^n}Z_0(h,x-y;x)[f(t-h,y)-f(t-h,x)]\,dy+f(t-h,x),
\end{multline}
due to (4.13). By Proposition 3, the absolute value of the first
integral in the right-hand side of (5.29) does not exceed
$$
Ch^{-\alpha}\int\limits_{\mathbb R^n}|x-y|^{-n+2+\gamma}\exp \left\{
-\sigma \left( h^{-\alpha
/2}|x-y|\right)^{\frac{2}{2-\alpha}}\right\}\,dy=\text{const}\cdot
h^{\frac{\gamma \alpha}2}\to 0,
$$
as $h\to 0$. Similarly, using (4.1) we show that the second
integral in (5.29) tends to zero as $h\to 0$. Thus
\begin{equation}
v_h^{(1)}(t,x)\longrightarrow f(t,x),
\end{equation}
as $h\to 0$.

Turning to $v_h^{(2)}$ we write
\begin{multline*}
v_h^{(2)}(t,x)=\int\limits_0^{t-h}d\lambda \int\limits_{\mathbb R^n}
\frac{\partial Z_0(t-\lambda ,x-y;y)}{\partial t}[f(\lambda ,y)-f(\lambda
,x)]\,dy\\
+\int\limits_0^{t-h}f(\lambda ,x)\,d\lambda \int\limits_{\mathbb R^n}
\frac{\partial Z_0(t-\lambda ,x-y;y)}{\partial t}\,dy.
\end{multline*}

It follows from (4.15) that
$$
\left| \int\limits_0^{t-h}d\lambda \int\limits_{\mathbb R^n}
\frac{\partial Z_0(t-\lambda ,x-y;y)}{\partial t}[f(\lambda ,y)-f(\lambda
,x)]\,dy\right| \le C(t-\lambda )^{-1+\frac{\alpha \gamma}2}.
$$
Together with (4.16) this shows that
\begin{multline}
v_h^{(2)}(t,x)\longrightarrow
\int\limits_0^td\lambda \int\limits_{\mathbb R^n}
\frac{\partial Z_0(t-\lambda ,x-y;y)}{\partial t}[f(\lambda ,y)-f(\lambda
,x)]\,dy\\
+\int\limits_0^tf(\lambda ,x)\,d\lambda \int\limits_{\mathbb R^n}
\frac{\partial Z_0(t-\lambda ,x-y;y)}{\partial t}\,dy.
\end{multline}

By the definition of $Y_0$, $\dfrac{\partial Z_0}{\partial t}=\D
Y_0$, and it follows from (5.30) and (5.31) that the fractional
derivative $\D W$ exists and can be represented as follows:
\begin{multline}
\left( \D W\right) (t,x)=f(t,x)+\int\limits_0^td\lambda
\int\limits_{\mathbb R^n}
\frac{\partial Z_0(t-\lambda ,x-y;y)}{\partial t}[f(\lambda ,y)-f(\lambda
,x)]\,dy\\
+\int\limits_0^tf(\lambda ,x)\,d\lambda \int\limits_{\mathbb R^n}
\frac{\partial Z_0(t-\lambda ,x-y;y)}{\partial t}\,dy.
\end{multline}

The case, in which $f(\lambda ,y)=Q(\lambda ,y;\xi )$, is
considered similarly, on the basis of the estimates $(4.15')$,
(4.16), (5.17), (5.18).

\medskip
{\bf 5.4.} {\it The initial condition.} The above study of heat
potentials, together with the investigation of the integral
$\int\limits_{\mathbb R^n}Z_0(t,x-\xi )u_0(\xi )\,d\xi $ performed
in \cite{K2} shows that our construction of the kernels $Z,Y$ indeed
gives, via the formula (1.4), a solution of the equation (1.1). It
remains to verify the initial condition (1.3).

By our construction,
$$
Z(t,x;\xi )=Z_0(t,x-\xi ,\xi )+V_Z(t,x;\xi ),
$$
where
$$
V_Z(t,x;\xi )=\int\limits_0^td\lambda \int\limits_{\mathbb R^n}
Y_0(t-\lambda ,x-y;y)Q(\lambda ,y;\xi )\,dy
$$
(see (5.1)). We have
\begin{multline*}
\int\limits_{\mathbb R^n}Z_0(t,x-\xi ,\xi )u_0(\xi )\,d\xi
=\int\limits_{\mathbb R^n}Z_0(t,x-\xi ,x)u_0(\xi )\,d\xi \\
+\int\limits_{\mathbb R^n}[Z_0(t,x-\xi ,\xi )-Z_0(t,x-\xi ,x)]u_0(\xi
)\,d\xi .
\end{multline*}
The first summand tends to $u_0(x)$, as $t\to 0$ \cite{K2}. By
Proposition 3, the second summand is majorized by
\begin{multline*}
Ct^{-\alpha }\int\limits_{\mathbb R^n}|x-\xi |^{-n+2+\gamma}\exp
\left\{ -\sigma \left( t^{-\alpha /2}|x-\xi
|\right)^{\frac{2}{2-\alpha}}\right\}\,d\xi \\
=Ct^{\gamma \alpha}\int\limits_{\mathbb R^n}|\eta |^{-n+2+\gamma}\exp
\left\{ -\sigma |\eta |^{\frac{2}{2-\alpha}}\right\}\,d\eta
\longrightarrow 0,
\end{multline*}
as $t\to 0$.

Let us find an estimate for $V_Z$. Suppose, for example, that
$n=3$ (other cases are treated similarly). Then by (4.7) and
(5.17)
$$
\left| V_Z(t,x;\xi )\right| \le C\int\limits_0^t(t-\lambda
)^{-\frac{\alpha}2-1}\lambda^{-\alpha}\,d\lambda \int\limits_{\mathbb
R^3}|y-\xi |^{-3+\gamma }\exp \{ -\sigma [\rho (t,x;\lambda
,y)+\rho (\lambda ,y;0,\xi )]\}\,dy.
$$
Taking $0<\gamma_0<\gamma$ we can rewrite this as
\begin{multline*}
\left| V_Z(t,x;\xi )\right| \le C\int\limits_0^t(t-\lambda
)^{-\frac{\gamma_0\alpha}2-1}\lambda^{-\alpha}\,d\lambda
\int\limits_{\mathbb R^3}\left( \frac{|x-y|}{(t-\lambda )^{\alpha
/2}}\right)^{1+\gamma_0}|x-y|^{-1-\gamma_0}
|y-\xi |^{-3+\gamma }\\
\times \exp \{ -\sigma [\rho (t,x;\lambda
,y)+\rho (\lambda ,y;0,\xi )]\}\,dy \\
\le C_1\int\limits_0^t(t-\lambda
)^{-\frac{\gamma_0\alpha}2-1}\lambda^{-\alpha}\,d\lambda
\int\limits_{\mathbb R^3}|x-y|^{-1-\gamma_0}
|y-\xi |^{-3+\gamma }\exp \{ -\sigma_1[\rho (t,x;\lambda
,y)+\rho (\lambda ,y;0,\xi )]\}\,dy,
\end{multline*}
and now Lemma 3 implies the estimate (iii) from the main Theorem.

Other estimates given in the formulation of our main Theorem
are proved in a quite similar way (note that we roughen some
estimates containing logarithmic terms replacing the logarithm
by the power function with an arbitrarily small exponent).
Returning now to the
verification of the initial condition we obtain that
\begin{multline*}
\left| \int\limits_{\mathbb R^3}V_Z(t,x;\xi )u_0(\xi )\,d\xi \right|
\le Ct^{\frac{\gamma_0\alpha}2-\alpha}\int\limits_{\mathbb
R^3}|\xi |^{\gamma -\gamma_0-1}\exp \left\{ -\sigma \left(
t^{-\alpha /2}|\xi |\right)^{\frac{2}{2-\alpha}}\right\}\,d\xi \\
=Ct^{\frac{\gamma\alpha}2}\int\limits_{\mathbb
R^3}|\eta |^{\gamma -\gamma_0-1}\exp \left\{ -\sigma
|\eta |^{\frac{2}{2-\alpha}}\right\}\,d\eta \longrightarrow 0,
\end{multline*}
as $t\to 0$. Thus, the initial condition (1.3) has been verified.

\medskip
{\bf 5.5.} {\it The H\"older continuity in $t$}. Let us prove the
assertion b) of the main theorem. Suppose that $0<\nu <1-\alpha$,
$\alpha +\lambda <\nu$. We have to show that for any fixed $x\in
\mathbb R^n$ $t^\nu u(t,x)$ is H\"older continuous in $t$ with the
exponent $\alpha +\lambda$. It is sufficient to prove the H\"older
continuity near the origin $t=0$. Let us consider, for example,
the case $n=3$; all other cases are treated similarly.

Let us write $u(t,x)=u_1(t,x)+u_2(t,x)+u_3(t,x)+u_4(t,x)$ where
\begin{gather*}
u_1(t,x)=\int\limits_{\mathbb R^3}Z_0(t,x-\xi ;\xi )u_0(\xi
)\,d\xi ,\\
u_2(t,x)=\int\limits_{\mathbb R^3}V_Z(t,x;\xi )u_0(\xi )\,d\xi ,\\
u_3(t,x)=\int\limits_0^td\lambda \int\limits_{\mathbb
R^3}Y_0(t-\lambda ,x-\xi ;\xi )f(\lambda ,\xi )\,d\xi ,\\
u_4(t,x)=\int\limits_0^td\lambda \int\limits_{\mathbb
R^3}V_Y(t-\lambda ,x;\xi )f(\lambda ,\xi )\,d\xi .
\end{gather*}

The H\"older continuity of $t^\nu u_1(t,x)$ was proved in
\cite{K2} (with the use of fractional calculus and the asymptotics
of H-functions) for the case of constant coefficients. The proof
for our case, in which there is also a dependence on the parameter
$\xi$, is identical to \cite{K2}.

As we saw above, $|u_2(t,x)|\le Ct^{\alpha \gamma /2}$, so that
$\left| t^\nu u_2(t,x)\right| \le Ct^{\frac{\alpha
\gamma}2+\nu}<Ct^{\alpha +\lambda}$ for small values of $t$, since
$\alpha +\lambda <\nu$.

Next, by (4.7)
$$
|u_3(t,x)|\le C\int\limits_0^t(t-\lambda )^{-\frac{\alpha}2-1}\,d\lambda
\int\limits_{\mathbb R^3}\exp \left\{ -\sigma \left( \frac{|x-\xi
|}{(t-\lambda )^{\alpha /2}}\right)^{\frac{2}{2-\alpha}}\right\}\,d\xi
\le C_1\int\limits_0^t(t-\lambda )^{\alpha -1}\,d\lambda
=C_2t^\alpha ,
$$
so that $\left| t^\nu u_3(t,x)\right| \le C_2t^{\alpha
+\nu}<C_2t^{\alpha +\lambda}$ for small values of $t$.

Finally, from the estimate (iii) of the main Theorem we find that
$$
\left| t^\nu u_4(t,x)\right| \le Ct^{\frac{\alpha \gamma}2+\alpha
+\nu}<Ct^{\alpha +\lambda}
$$
for small values of $t$, as desired.

\medskip
\section{THE ONE-DIMENSIONAL CASE}

In the case $n=1$ the parametrix and the fundamental solution have
no singularity in the spatial variable. This simplifies the
situation greatly and makes it similar to the theory of conventional
parabolic equations of an arbitrary order. Note however that even
this case is more complicated than the classical study of a
second order parabolic differential equation. In particular, the
estimates of iterated kernels in the Levi method should still be
performed in two stages (see Lemma 6).

Here we give the main estimates for this case; their proofs,
simplified versions of those given the preceding section, are
omitted. The scheme of the Levi method and the main notations are
as above; see (5.1)-(5.4).

\begin{prop}
If $n=1$, then
\begin{gather*}
|M(t,x;\xi )|\le Ct^{-(3-\gamma )\alpha /2}\exp \{-\sigma \rho
(t,x;0,\xi )\},\\
|K(t,x;\xi )|\le Ct^{-1-(1-\gamma )\alpha /2}\exp \{-\sigma \rho
(t,x;0,\xi )\},\\
|Q(t,x;\xi )|\le Ct^{-(3-\gamma )\alpha /2}\exp \{-\sigma \rho
(t,x;0,\xi )\},\\
|\Psi (t,x;\xi )|\le Ct^{-1-(1-\gamma )\alpha /2}\exp \{-\sigma \rho
(t,x;0,\xi )\}.
\end{gather*}
\end{prop}

\medskip
The estimates for differences of the above kernels are collected
in the next proposition.

\begin{prop}
If $n=1$, then
\begin{gather*}
|\Delta_xM(t,x;\xi )|\le C|x-x'|^{\gamma -\varepsilon}
t^{-(3-\varepsilon )\alpha /2}\exp \{-\sigma \rho
(t,x'';0,\xi )\},\\
|\Delta_xK(t,x;\xi )|\le C|x-x'|^{\gamma -\varepsilon}
t^{-1-(1-\varepsilon )\alpha /2}\exp \{-\sigma \rho
(t,x'';0,\xi )\},\\
|\Delta_xQ(t,x;\xi )|\le C|x-x'|^{\gamma -\varepsilon}
t^{-(3-\varepsilon )\alpha /2}\exp \{-\sigma \rho
(t,x'';0,\xi )\},\\
|\Delta_x\Psi (t,x;\xi )|\le C|x-x'|^{\gamma -\varepsilon}
t^{-1-(1-\varepsilon )\alpha /2}\exp \{-\sigma \rho
(t,x'';0,\xi )\},
\end{gather*}
$0<\varepsilon <\gamma$.
\end{prop}

\medskip
\section{NONNEGATIVITY}

Let us prove that the functions $Z$ and $Y$ are nonnegative. It is
sufficient to show the nonnegativity of the solution (1.4) of the
problem (1.1), (1.3) for arbitrary nonnegative functions $u_0\in
C_0^\infty (\mathbb R^n)$, $f\in C_0^\infty ([0,T]\times \mathbb
R^n)$.

As we already know, the solution $u(t,x)$ is bounded (due to our
estimates of $Z$ and $Y$), say $|u(t,x)|\le M$, and belongs to
$H_\nu ^{\alpha +\lambda }[0,T]$, $0<\lambda <1-\alpha$, $\alpha
+\lambda <\nu$, for each $x$. This regularity in $t$ makes it
possible to represent the regularized fractional derivative $\D u$
via the Marchaud fractional derivative:
$$
\left( \D u\right) (t,x)=\frac{1}{\Gamma (1-\alpha )}\left[
t^{-\alpha}u(t,x)-t^{-\alpha}u_0(x)\right]
+\frac{\alpha}{\Gamma (1-\alpha )}\lim\limits_{\varepsilon \to
+0}\psi_\varepsilon (t,x)
$$
where $\psi_\varepsilon (t,x)=0$ for $0<t\le \varepsilon$,
$$
\psi_\varepsilon
(t,x)=\int\limits_0^{t-\varepsilon}\frac{u(t,x)-u(\tau
,x)}{(t-\tau )^{1+\alpha }}\,d\tau ,\quad \varepsilon <t\le T,
$$
and the limit exists for all $x\in \mathbb R^n$, $t\in [0,T]$ (see
Sect. 13 in \cite{SKM}).

Denote $v(t,x)=\dfrac{u(t,x)}{E_\alpha (\beta t^\alpha )}$ where
$E_\alpha$ is the Mittag-Leffler function \cite{Er}, $\beta >0$
will be chosen later. Using the fact that $\D E_\alpha (\beta t^\alpha
)=\beta E_\alpha (\beta t^\alpha )$ we obtain that
$$
\D u(t,x)=\beta E_\alpha (\beta t^\alpha
)v(t,x)-\frac{t^{-\alpha}}{\Gamma (1-\alpha )}u_0(x)+\left(
L_\beta v\right) (t,x)
$$
where
$$
\left( L_\beta v\right) (t,x)=\frac{1}{\Gamma (1-\alpha )}\left[
v(t,x)t^{-\alpha}+\alpha \lim\limits_{\varepsilon \to +0}
\int\limits_0^{t-\varepsilon}E_\alpha (\beta \tau^\alpha
)\frac{v(t,x)-v(\tau ,x)}{(t-\tau )^{1+\alpha }}\,d\tau \right].
$$

Now the equation (1.1) takes the form
\begin{equation}
\frac{1}{E_\alpha (\beta t^\alpha )}\left( L_\beta v\right)
(t,x)-(B-\beta )v(t,x)=g(t,x)
\end{equation}
where
$$
g(t,x)=\frac{1}{E_\alpha (\beta t^\alpha )}\left[
\frac{t^{-\alpha}}{\Gamma (1-\alpha )}u_0(x)+f(t,x)\right].
$$

Let $c_0=\sup\limits_{x\in \mathbb R^n}|c(x)|$, $d>0$, $\beta
=c_0+d$. Consider the function
$$
F_R(t,x)=\frac{M}{R^2}\left( |x|^2+\mu t^\alpha +1\right) ,\quad
R,\mu >0.
$$
Since $F_R(t,x)$ and $E_\alpha (\beta t^\alpha )$ are monotone
increasing in $t$, we have
$$
\left( L_\beta F_R\right) (t,x)\ge \left( D_{0+}^\alpha F_R\right)
(t,x)\ge \frac{M\mu }{R^2}D_{0+}^\alpha (t^\alpha )=
\frac{M\mu \Gamma (1+\alpha )}{R^2}.
$$
Here $D_{0+}^\alpha$ is the Riemann-Liouville fractional
derivative (see \cite{SKM}), and in estimating $L_\beta F_R$ we
used the Marchaud form of $D_{0+}^\alpha$.

On the other hand,
\begin{multline*}
(B-\beta )F_R(t,x)=\sum\limits_{i=1}^na_{ii}(x)\frac{2M}{R^2}
+\sum\limits_{j=1}^nb_j(x)\frac{2M}{R^2}x_j+[c(x)-c_0]F_R(t,x)-dF_R(t,x)\\
\le \frac{M}{R^2}\left( C_1+C_2|x|-d|x|^2-d\mu t^\alpha -d\right),
\end{multline*}
so that
\begin{equation}
\frac{1}{E_\alpha (\beta t^\alpha )}\left( L_\beta F_R\right)
(t,x)-(B-\beta )F_R(t,x)\ge
\frac{M}{R^2}\left[ \frac{\mu \Gamma (1+\alpha )}{E_\alpha (\beta t^\alpha )}
-C_1-C_2|x|+d|x|^2+d\mu t^\alpha +d\right] \ge 0
\end{equation}
for all $x\in \mathbb R^n$, $t\in [0,T]$, if the number $\mu$ is
large enough.

Let $G(t,x)=v(t,x)+F_R(t,x)$. Since $g(t,x)\ge 0$, it follows from
(7.1) and (7.2) that
\begin{equation}
\frac{1}{E_\alpha (\beta t^\alpha )}\left( L_\beta G\right)
(t,x)-(B-\beta )G(t,x)\ge 0.
\end{equation}
If $|x|=R$, then
$$
G(t,x)=v(t,x)+M+\frac{\mu Mt^\alpha }{R^2}+\frac{M}{R^2}>0.
$$
Also $G(0,x)=u_0(x)+F_R(0,x)>0$. Then $G(t,x)\ge 0$ for all
$(t,x)$ with $t\in [0,T]$, $|x|<R$. Indeed, otherwise the function
$G$ possesses a point $(t^0,x^0)$ of its global minimum on the set
$\{ 0<t\le T,|x|<R\}$, such that $G(t^0,x^0)<0$. At this point
$(B-\beta )G(t^0,x^0)\ge 0$ \cite{F}, and by (7.3) $\left( L_\beta G\right)
(t^0,x^0)\ge 0$, which contradicts the definition of $L_\beta$.

Thus, if $|x|\le R$, then
$$
u(t,x)\ge -\frac{M}{R^2}\left( |x|^2+\mu t^\alpha +1\right)
E_\alpha (\beta t^\alpha ) .
$$
Since $R$ is arbitrary, this means that $u(t,x)\ge 0$.


\medskip
\begin{thebibliography}{999}
\bibitem{AL}
V. V. Anh and N. N. Leonenko, Spectral analysis of fractional
kinetic equations with random data, {\it J. Statist. Phys.} {\bf
104} (2001), 1349--1387.
\bibitem{BM}
B. Baeumer and M. Meerschaert, Stochastic solutions for fractional
Cauchy problems, {\it Fract. Calc. Appl. Anal.} {\bf 4}
(2001), 481--500.
\bibitem{Bazh}
E. Bazhlekova, The abstract Cauchy problem for fractional
evolution equation, {\it Fract. Calc. Appl. Anal.} {\bf 1}
(1998), 255--270.
\bibitem{Bazh1}
E. Bazhlekova, Fractional Evolution Equations in Banach Spaces,
Dissertation, Technische Universiteit Eindhoven, 2001.
\bibitem{Br}
B. L. J. Braaksma, Asymptotic expansions and analytic
continuation for a class of Barnes integrals, {\it Compositio
Math.} {\bf 15} (1964), 239--341.
\bibitem{DN}
M. M. Dzhrbashyan and A. B. Nersessyan, Fractional derivatives
and Cauchy problem for differential equations of fractional
order, {\it Izv. AN Arm. SSR. Matematika} {\bf 3} (1968), 3--29
(Russian).
\bibitem{E}
S. D. Eidelman, Parabolic Systems, North-Holland,
Amsterdam, 1969.
\bibitem{El}
A. M. El-Sayed, Fractional order evolution equations, {\it J.
Fract. Calc.} {\bf 7} (1995), 89--100.
\bibitem{Er}
A. Erdelyi, W. Magnus, F. Oberhettinger, and F. Tricomi,
Higher Transcendental Functions. Vol. III, McGraw-Hill, New York,
1955.
\bibitem{F}
A. Friedman, Partial Differential Equations of Parabolic
Type, Prentice-Hall, Englewood Cliffs, NJ, 1964.
\bibitem{Gor}
R. Gorenflo, F. Mainardi, D. Moretti, and P. Paradisi, Time fractional
diffusion: A discrete random walk approach, {\it Nonlinear Dynamics}
{\bf 29} (2002), 129--143.
\bibitem{K1}
A. N. Kochubei, A Cauchy problem for evolution equations of
fractional order, {\it Differential Equations} {\bf 25} (1989),
967--974.
\bibitem{K2}
A. N. Kochubei, Fractional-order diffusion, {\it Differential Equations}
{\bf 26} (1990), 485--492.
\bibitem{Kols}
T. Kolsrud, On a class of probabilistic integrodifferential
equations. In: {\it Ideas and Methods in Mathematics and Physics.
Memorial Volume Dedicated to Raphael H\o egh-Krohn, Vol. 1},
Cambridge University Press, 1992, pp. 168--172.
\bibitem{Kost}
V. A. Kostin, Cauchy problem for an abstract differential
equation with fractional derivatives, {\it Russian Acad. Sci.
Dokl. Math.} {\bf 46} (1993), 316--319.
\bibitem{LSU}
O. A. Ladyzhenskaya, V. A. Solonnikov, and N. N. Uraltseva,
Linear and Quasilinear Equations of Parabolic Type, American
Mathematical Society, Providence, 1968.
\bibitem{Meer}
M. M. Meerschaert, D. A. Benson, H. P. Scheffler and B. Baeumer,
Stochastic solutions of space-time fractional diffusion equations,
{\it Phys. Rev. E} {\bf 65} (2002), 1103--1106.
\bibitem{MK}
R. Metzler and J. Klafter, The random walk's guide to anomalous
diffusion: a fractional dynamics approach, {\it Physics Reports},
{\bf 339} (2000), 1--77.
\bibitem{MR}
K. Miller and B. Ross, An Introduction to the Fractional Calculus
and Fractional Differential Equations, Wiley and Sons, New York,
1993.
\bibitem{PBM}
A. P. Prudnikov, Yu. A. Brychkov, and O. I. Marichev, Integrals
and Series. Vol. 3: More Special Functions, Gordon and Breach,
New York, 1990.
\bibitem{SKM}
S. G. Samko, A. A. Kilbas, and O. I. Marichev,
Fractional Integrals and Derivatives: Theory and Applications,
Gordon and Breach, New York, 1993.
\bibitem{SW}
W. R. Schneider, W. Wyss, Fractional diffusion and wave
equations, {\it J. Math. Phys.} {\bf 30} (1989), 134--144.
\bibitem{S1}
W. R. Schneider, Fractional diffusion, {\it Lecture Notes Phys.},
{\bf 355} (1990), 276--286.
\bibitem{S2}
W. R. Schneider, Grey noise. In: {\it Ideas and Methods in Mathematics and
Physics. Memorial Volume Dedicated to Raphael H\o egh-Krohn, Vol. 1},
Cambridge University Press, 1992, pp. 261--282.
\bibitem{SGG}
H. M. Srivastava, K. C. Gupta, and S. P. Goyal, The H-Functions
of One and Two variables with Applications, South Asian
Publishers, New Dehli, 1982.
\bibitem{W}
W. Wyss, The fractional diffusion equation, {\it J. Math. Phys.}
{\bf 27} (1986), 2782--2785.
\bibitem{Yor}
M. Yor, W. Schneider's grey noise and fractional Brownian motion.
In: {\it Proc. Easter Meeting on Probability (Edinburgh, April
10-14, 1989)}.
\end{thebibliography}
\end{document}